\documentclass[11pt]{article}
\usepackage{mathrsfs}
\usepackage{}
\usepackage[hang,flushmargin]{footmisc}
\usepackage{multirow}
\usepackage{amssymb}
\usepackage{amsfonts}
\usepackage{caption}
\usepackage{amsmath}
\usepackage{graphicx}
\usepackage{amsmath,amssymb,latexsym,color}
\usepackage[mathscr]{eucal}
 \makeatletter
    
    \newcommand{\Rmnum}[1]
    {\expandafter\@slowromancap\romannumeral #1@}
    \makeatother
\def\wz{\tilde}

\newtheorem{thm}{Theorem}[section]

\newtheorem{prop}[thm]{Proposition}
\newtheorem{lemma}[thm]{Lemma}
\newcounter{foo}[subsection]
\newcounter{fooo}[section]
\newtheorem{step}[foo]{Step}
\newtheorem{stepp}[fooo]{Step}

\newtheorem{example}[thm]{Example}
\newtheorem{defin}[thm]{Definition}

\newtheorem{remark}{Remark}[section]

\newcommand{\qed}{\hfill\Box\medskip}

\usepackage{CJK}
\begin{document}
\begin{CJK*}{GBK}{song}

\renewcommand{\baselinestretch}{1.3}
\title{Quasi-thin weakly distance-regular digraphs}

\author{
Yuefeng Yang\textsuperscript{a,b}\quad Benjian Lv\textsuperscript{b}\quad
Kaishun Wang\textsuperscript{b}\\
\\{\footnotesize  \textsuperscript{a} \em  School of Science, China University of Geosciences, Beijing, 100083, China}\\{\footnotesize  \textsuperscript{b} \em Sch. Math. Sci. {\rm \&} Lab. Math. Com. Sys., Beijing Normal University, Beijing, 100875, China }  }
\date{}
\maketitle
\footnote{\scriptsize
{\em E-mail address:} yangyf@mail.bnu.edu.cn(Y.Yang), bjlv@bnu.edu.cn(B.Lv), wangks@bnu.edu.cn(K.Wang).}

\begin{abstract}
A weakly distance-regular digraph is quasi-thin if the maximum value of its intersection numbers is $2$. In this paper, we focus on commutative quasi-thin weakly distance-regular
digraphs, and classify such digraphs with valency more than $3$. As a result, this family of
digraphs are completely determined.

\medskip
\noindent {\em AMS classification:} 05E30

\noindent {\em Key words:} Weakly distance-regular digraph;  quasi-thin; Cayley digraph

\end{abstract}

\section{Introduction}

Throughout this paper $\Gamma$ always denotes a finite simple digraph. We write $V\Gamma$ and $A\Gamma$ for the vertex set and arc set of $\Gamma$, respectively. A \emph{path} of length $r$ from $x$ to $y$ is a sequence of vertices $(x=w_{0},w_{1},\ldots,w_{r}=y)$
such that $(w_{t-1}, w_{t})\in A\Gamma$ for $t=1,2,\ldots,r$. A digraph is said to be \emph{strongly connected} if, for any two distinct vertices $x$ and $y$, there is a path from $x$ to $y$. The length of a shortest path from $x$ to $y$ is called the \emph{distance} from $x$ to $y$
in $\Gamma$, denoted by $\partial_\Gamma(x,y)$. Let $\wz{\partial}_\Gamma(x,y)=(\partial_\Gamma(x,y),\partial_\Gamma(y,x))$ and $\wz{\partial}(\Gamma)=\{\wz{\partial}_\Gamma(x,y)\mid x,y\in V\Gamma\}$. We call $\wz{\partial}_\Gamma(x,y)$ the \emph{two way distance} from $x$ to $y$ in $\Gamma$.
If no confusion occurs, we  write $\partial(x,y)$ (resp. $\wz{\partial}(x,y)$) instead of $\partial_\Gamma(x,y)$ (resp. $\wz{\partial}_\Gamma(x,y)$).
An arc $(u,v)$ of $\Gamma$ is of \emph{type} $(1,r)$ if $\partial(v,u)=r$. A path $(w_{0},w_{1},\ldots,w_{r-1})$ is said to be a \emph{circuit} of length
$r$ if
$\partial(w_{r-1},w_{0})=1$. A circuit is \emph{undirected} if each of its arcs is of type $(1,1)$. Let $C_{r}$ denote the undirected circuit of length $r$.

A strongly
connected digraph $\Gamma$ is said to be \emph{weakly distance-regular} if, for any $\wz{h}$, $\wz{i}$, $\wz{j}\in\wz{\partial}(\Gamma)$, the cardinality of the set $$P_{\wz{i},\wz{j}}(x,y):=\{z\in V\Gamma\mid\wz{\partial}(x,z)=\wz{i}~\textrm{and}~\wz{\partial}(z,y)=\wz{j}\}$$ is constant whenever $\wz{\partial}(x,y)=\wz{h}$. This constant is denoted by $p_{\wz{i},\wz{j}}^{\wz{h}}$. The integers
$p_{\wz{i},\wz{j}}^{\wz{h}}$ are called the \emph{intersection numbers}. We say
that $\Gamma$ is \emph{commutative} if $p_{\wz{i},\wz{j}}^{\wz{h}}=p_{\wz{j},\wz{i}}^{\wz{h}}$ for all $\wz{i}$, $\wz{j}$,
$\wz{h}\in\wz{\partial}(\Gamma)$. A weakly distance-regular digraph is \emph{quasi-thin} (resp. \emph{thin}) if the maximum value of its intersection
numbers is $2$ (resp. $1$). The size of $\Gamma_{\tilde{i}}(x):=\{y\in V\Gamma\mid\tilde{\partial}(x,y)=\tilde{i}\}$
depends only on $\tilde{i}$, denoted by $k_{\tilde{i}}$. The integer $k:=\sum_{(1,j)\in\wz{\partial}(\Gamma)}k_{1,j}$ is called the \emph{valency} of $\Gamma$, which is often called the \emph{out-degree} of $\Gamma$.

Some special families of weakly distance-regular digraphs were classified. See \cite{KSW03,HS04} for valency $2$, \cite{KSW04,YYF16,YYF} for valency $3$ and \cite{HS04} for thin case.
In this paper, we classify commutative quasi-thin weakly distance-regular digraphs of valency more than $3$, and obtain the following main result.

\begin{thm}\label{Main}
If $\Gamma$ is a
commutative quasi-thin weakly distance-regular digraph of valency more than $3$, then $\Gamma$ is isomorphic to one of the following Cayley
digraphs:\vspace{-0.3cm}
\begin{itemize}
\item [{\rm(i)}] ${\rm Cay}(\mathbb{Z}_{8},\{1,2,3,6\})$.\vspace{-0.3cm}

\item [{\rm(ii)}] ${\rm Cay}(\mathbb{Z}_{4p},\{1,2,2p+i,2p+1,2p+2\})$, $p\neq 2-i$.\vspace{-0.3cm}

\item [{\rm(iii)}] ${\rm Cay}(\mathbb{Z}_{4}\times\mathbb{Z}_{4},\{(0,1),(1,0),(2,0),(0,2)\})$.\vspace{-0.3cm}

\item [{\rm(iv)}] ${\rm Cay}(\mathbb{Z}_{q}\times\mathbb{Z}_{4},\{(0,1),(1,0),(1,2),(0,2+i)\})$, $q\neq3+i$.\vspace{-0.3cm}

\item [{\rm(v)}] ${\rm Cay}(\mathbb{Z}_{2q}\times\mathbb{Z}_{2},\{(0,1),(1,0),(2,0),(1,1)\})$.\vspace{-0.3cm}

\item [{\rm(vi)}] ${\rm Cay}(\mathbb{Z}_{4q}\times\mathbb{Z}_{2},\{(0,1),(1,0),(2,0),(2q+1,0),(2q+2,0),(2qi,1)\})$, $q\notin\{3,3+i\}$.\vspace{-0.3cm}

\item [{\rm(vii)}] ${\rm Cay}(\mathbb{Z}_{2q}\times\mathbb{Z}_{4},\{(0,1),(1,0),(1,2),(0,2-i),(2,0),(2,2)\})$, $q\notin\{3,3+i\}$.\vspace{-0.3cm}

\item [{\rm(viii)}] ${\rm Cay}(\mathbb{Z}_{2q}\times\mathbb{Z}_{n},\{(0,1),(1,0),(2,0),(0,-1)\})$.\vspace{-0.3cm}

\item [{\rm(ix)}] ${\rm Cay}(\mathbb{Z}_{2q}\times\mathbb{Z}_{n},\{(0,1),(1,(c+1)/2),(1,(c-1)/2),(2,c),(0,-1)\})$.\vspace{-0.3cm}

\item [{\rm(x)}] ${\rm Cay}(\mathbb{Z}_{2n}\times\mathbb{Z}_{q},\{(0,1),(1,(t+1)/2),(-1,(1-t)/2),(2,t),(-2,-t)\})$.\vspace{-0.3cm}
\end{itemize}
Here, $i\in\{0,1\}$, $2\leq p$, $3\leq q$, $3\leq n\leq q-(1+(-1)^{q})/2$, $c=n/{\rm gcd}(q,n)$, $t=q/{\rm gcd}(q,n)$ and $c,t$ are both odd.
\end{thm}

Routinely, all digraphs in above theorem are commutative quasi-thin weakly distance-regular. For the last seven families of Cayley digraphs, in Table 1, we list the two way distance from the identity element to any other element of the corresponding group.

In order to give a high level description of our proof of Theorem \ref{Main}, we need additional notations and terminologies. Let $\Gamma$ be a weakly distance-regular digraph and $R=\{\Gamma_{\wz{i}}\mid\wz{i}\in\wz{\partial}(\Gamma)\}$, where $\Gamma_{\wz{i}}=\{(x,y)\in V\Gamma\times V\Gamma\mid\wz{\partial}(x,y)=\wz{i}\}$. Then $(V\Gamma,R)$ is an association scheme (\cite{EB84,PHZ96,PHZ05}). Moreover, if $\Gamma$ is quasi-thin, then $(V\Gamma,R)$ is quasi-thin. About this special scheme, see
 \cite{MH01,MH02,MM12}. For two nonempty subsets $E$ and $F$ of $R$, define
\begin{eqnarray*}
EF&:=&\{\Gamma_{\wz{h}}\mid\sum_{\Gamma_{\wz{i}}\in E}\sum_{\Gamma_{\wz{j}}\in F}p_{\wz{i},\wz{j}}^{\wz{h}}\neq0\},
\end{eqnarray*}
and write $\Gamma_{\wz{i}}\Gamma_{\wz{j}}$ instead of $\{\Gamma_{\wz{i}}\}\{\Gamma_{\wz{j}}\}$. For any $(a,b)\in\wz{\partial}(\Gamma)$, we usually write $k_{a,b}$ (resp. $\Gamma_{a,b}$) instead of $k_{(a,b)}$ (resp. $\Gamma_{(a,b)}$). Now we list basic properties of intersection numbers which are used frequently in this paper.

\begin{lemma}\label{jiben}
{\rm (\cite[Chapter \Rmnum{2}, Proposition 2.2]{EB84} and \cite[Proposition 5.1]{ZA99})} For each $\wz{i}:=(a,b)\in\wz{\partial}(\Gamma)$, define $\wz{i}^{*}=(b,a)$. The following hold:\vspace{-0.3cm}
\begin{itemize}
\item [{\rm(i)}] $k_{\wz{d}}\;k_{\wz{e}}=\sum_{\wz{f}\in\wz{\partial}(\Gamma)}p_{\wz{d},\wz{e}}^{\wz{f}}\;k_{\wz{f}}$.\vspace{-0.3cm}

\item [{\rm(ii)}] $p_{\wz{d},\wz{e}}^{\wz{f}}\;k_{\wz{f}}=p_{\wz{f},\wz{e}^{*}}^{\wz{d}}\;k_{\wz{d}}=p_{\wz{d}^{*},\wz{f}}^{\wz{e}}\;k_{\wz{e}}$.\vspace{-0.3cm}

\item [{\rm(iii)}] $|\Gamma_{\wz{d}}\;\Gamma_{\wz{e}}|\leq{\rm gcd}(k_{\wz{d}},k_{\wz{e}})$.\vspace{-0.3cm}

\item [{\rm(iv)}] $\sum_{\wz{e}\in\wz{\partial}(\Gamma)}p_{\wz{d},\wz{e}}^{\wz{f}}=k_{\wz{d}}$.\vspace{-0.3cm}

\item [{\rm(v)}] ${\rm lcm}(k_{\wz{d}},k_{\wz{e}})\mid p_{\wz{d},\wz{e}}^{\wz{f}}\;k_{\wz{f}}$.\vspace{-0.3cm}

\item [{\rm(vi)}] $\sum_{\wz{f}\in\wz{\partial}(\Gamma)}p_{\wz{d},\wz{e}}^{\wz{f}}\;p_{\wz{g},\wz{f}}^{\wz{h}}=\sum_{\wz{l}\in\wz{\partial}(\Gamma)}p_{\wz{g},\wz{d}}^{\wz{l}}\;p_{\wz{l},\wz{e}}^{\wz{h}}$.
\end{itemize}
\end{lemma}

We now introduce the concepts about arcs. An arc of type $(1,q-1)$ is said to be \emph{pure}, if every circuit of length $q$ containing it consists of arcs of type $(1,q-1)$; otherwise, this arc is said to be \emph{mixed}. We say that $(1,q-1)$ is pure if any arc of type $(1,q-1)$ is pure; otherwise, we say that $(1,q-1)$ is mixed. The concepts of pure arc and mixed arc are inspired by Suzuki in \cite{HS04}.

Another concept we need is a configuration. Let $h$ and $q$ be distinct integers more than $2$. If $(\Gamma_{1,q-1})^{2}=\{\Gamma_{2,q-2}\}$ and $(\Gamma_{1,h-1})^{2}\subseteq\Gamma_{1,q-1}\Gamma_{q-1,1}$, we say that the configuration $C_{q,h}$ exists.

For fixed $x\in V\Gamma$, let $\Delta_{q_{1},q_{2},\ldots,q_{l}}(x)$ be the connected component of digraph $(V\Gamma,\cup_{i=1}^{l}\Gamma_{1,q_{i}-1})$ containing vertex $x$. Note that $\Delta_{q_{1},q_{2},\ldots,q_{l}}(x)$ does not depend on the choice of vertex $x$ up to isomorphism. If no confusion occurs, we write $\Delta_{q_{1},q_{2},\ldots,q_{l}}$ instead of $\Delta_{q_{1},q_{2},\ldots,q_{l}}(x)$.

Let $\Gamma$ be a commutative quasi-thin weakly distance-regular digraph of valency more than $3$ in the remaining of this paper. We are now ready to give a high level description of our proof of Theorem \ref{Main}.

\vspace{3ex}

\noindent\textbf{Outline of the proof of Theorem \ref{Main}.}

In Section 2, we give a characterization of mixed arcs, i.e. we show that $(1,q-1)$ is mixed if and only if $p_{(1,q-1),(1,q-1)}^{(1,q-2)}\neq0$ and $(1,q-2)$ is pure.

In Section 3, we discuss the basic properties about the configuration $C_{q,h}$. In particular, we show that, if $C_{q,h}$ exists, then $(1,q-1)$ is pure, $h$ is a constant and $h\in\{3,4\}$.

In Section 4, applying the results in Sections 2 and 3, we prove the following result.

\begin{prop}\label{jiegou}
Let $K=\{(1,r)\mid(1,r)\in\wz{\partial}(\Gamma)\}$. Then one of the following holds:\vspace{-0.3cm}
\begin{itemize}
\item [{\rm C1)}] $K=\{(1,1),(1,2),(1,q-1)\}$, where $C_{q,3}$ exists.\vspace{-0.3cm}

\item [{\rm C2)}] $K=\{(1,3),(1,q-1),(1,q)\}$, where $C_{q,4}$ exists and $(1,q)$ is mixed.\vspace{-0.3cm}

\item [{\rm C3)}] $K=\{(1,1),(1,2),(1,q-1),(1,q)\}$, where $C_{q,3}$ exists and $(1,q)$ is mixed.\vspace{-0.3cm}

\item [{\rm C4)}] $K=\{(1,1),(1,q-1)\}$, where $(1,q-1)$ is pure.\vspace{-0.3cm}

\item [{\rm C5)}] $K=\{(1,q-1),(1,q)\}$, where $(1,q)$ is mixed.\vspace{-0.3cm}

\item [{\rm C6)}] $K=\{(1,1),(1,q-1),(1,q)\}$, where $(1,q)$ is mixed.
\end{itemize}
\end{prop}
\vspace{-1ex}

In Section 5, we determine the subdigraphs $\Delta_{q,3}$ for the cases C1 and C3, the subdigraphs $\Delta_{q,4}$ for case C2, the subdigraphs $\Delta_{2,q}$ for cases C4 and C6, and the subdigraphs $\Delta_{q,q+1}$ for cases C5 and C6.

In Section 6, we give a proof of Theorem \ref{Main}. For the cases C1, C2 and C3, we determine $\Gamma$ based on the subdigraphs $\Delta_{q,3}$ and $\Delta_{q,4}$. For the cases C4, C5 and C6, we determine $\Gamma$ based on the subdigraphs $\Delta_{2,q}$ and $\Delta_{q,q+1}$.

\section{Characterization of mixed arcs}

The main result of this section is the following important result which characterizes mixed arcs.

\begin{thm}\label{Main1}
Let $q\geq3$ and $(1,q-1)\in\wz{\partial}(\Gamma)$.\vspace{-0.3cm}
\begin{itemize}
\item [{\rm(i)}] If $p_{(1,s-1),(1,t-1)}^{(1,q-1)}\neq0$ with $s<t$, then $s=2$ and $t=q$.\vspace{-0.3cm}

\item [{\rm(ii)}] The following are equivalent:\vspace{-0.3cm}

{\rm(a)} $(1,q-1)$ is mixed; {\rm(b)} $p_{(1,q-1),(1,q-1)}^{(1,q-2)}\neq0$ and $(1,q-2)$ is pure; {\rm(c)} $p_{(1,q-1),(1,q-1)}^{(1,s-1)}\neq0$ for some $s$.\vspace{-0.3cm}

\item [{\rm(iii)}] If $p_{(1,q-1),(1,q-1)}^{(1,s-1)}\neq0$, then $s=q-1$.
\end{itemize}
\end{thm}

In the proof of Theorem \ref{Main1}, we use the following auxiliary lemmas.

\begin{lemma}\label{jiben2}
Suppose $\wz{d},\wz{h},\wz{l}\in\wz{\partial}(\Gamma)$ and $k_{\wz{d}}=2$. The following hold:\vspace{-0.3cm}
\begin{itemize}
\item [{\rm(i)}] $k_{\wz{h}}=k_{\wz{h}^{*}}\leq2$.\vspace{-0.3cm}

\item [{\rm(ii)}] $|\Gamma_{\wz{h}}\Gamma_{\wz{l}}|\leq2$ and equality holds only if $k_{\wz{h}}=k_{\wz{l}}=2$.\vspace{-0.3cm}

\item [{\rm(iii)}] $p_{\wz{d},\wz{d}}^{\wz{e}}=2$ for some $\wz{e}\in\wz{\partial}(\Gamma)$.\vspace{-0.3cm}

\item [{\rm(iv)}] $\Gamma_{\wz{d}}\Gamma_{\wz{d}^{*}}=\{\Gamma_{0,0},\Gamma_{e,e}\}$. In particular, if $p_{\wz{d},\wz{d}^{*}}^{\wz{e}}\neq0$, then $\wz{e}=\wz{e}^{*}.$
\end{itemize}
\end{lemma}
\textbf{Proof.}~Since $k_{\wz{h}^{*}}=k_{\wz{h}}=p_{\wz{h},\wz{h}^{*}}^{(0,0)}$ by Lemma \ref{jiben} (ii), (i) is valid. (ii) follows from (i) and Lemma \ref{jiben} (iii). By the commutativity of $\Gamma$, (iii) holds. In view of (ii) and Lemma \ref{jiben} (i), (iv) is valid.$\qed$

The commutativity of $\Gamma$ will be used frequently in the sequel, so we no longer refer to it for the sake of simplicity.

\begin{lemma}\label{pure-g=2}
If $(x_{0},x_{1},\ldots,x_{n-1})$ is an undirected circuit in $\Gamma$, then $\partial(x_{0},x_{i})=\partial(x_{i},x_{0})=\partial(x_{0},x_{n-i})$ for $1\leq i\leq n-1$.
\end{lemma}
\textbf{Proof.}~It is routine by induction.$\qed$

\begin{lemma}\label{jiben3}
Let $q\geq3$. Suppose that $(1,q-1)$ is pure and $k_{1,q-1}=2$. Then one of the following holds:\vspace{-0.3cm}
\begin{itemize}
\item [{\rm(i)}] $p_{(1,q-1),(1,q-1)}^{(2,q-2)}=2$, $\Delta_{q}\simeq{\rm Cay}(\mathbb{Z}_{2q},\{1,q+1\})$, $(\Gamma_{1,q-1})^{i}=\{\Gamma_{i,q-i}\}$ for $2\leq i\leq q-1$.\vspace{-0.3cm}

\item [{\rm(ii)}] $p_{(1,q-1),(1,q-1)}^{(2,q-2)}=1$, $\Delta_{q}\simeq{\rm Cay}(\mathbb{Z}_{q}\times\mathbb{Z}_{q},\{(1,0),(0,1)\})$, $|(\Gamma_{1,q-1})^{2}|=2$.
\end{itemize}
\end{lemma}
\textbf{Proof}. Similar to the proofs of Lemma 12 in \cite{YYF16} and Proposition 4.3 in \cite{KSW03}.$\qed$

\begin{lemma}\label{mix-(2,q-2)'}
Let $q\geq3$. Suppose that $p_{(1,q-1),(1,q-1)}^{(1,q-2)}\neq0$ and $(1,q-2)$ is pure. Then the following hold:\vspace{-0.3cm}
\begin{itemize}
\item [{\rm(i)}] $p_{(1,q-1),(1,q-2)}^{(2,q-2)}\neq0$ and $p_{(1,q-1),(1,q-1)}^{(2,q-2)}=0$.\vspace{-0.3cm}

\item [{\rm(ii)}] Any circuit of length $q$ containing an arc of type $(1,q-1)$ consists of arcs of types $(1,q-1)$ and $(1,q-2)$. \vspace{-0.3cm}

\item [{\rm(iii)}] If $|(\Gamma_{1,q-1})^{2}|=2$ and $k_{1,q-2}=1$, then $p_{(1,q-1),(1,q-1)}^{(2,q-1)}\neq0$.
\end{itemize}
\end{lemma}
\textbf{Proof.}~(i)~Let $(z,z_{0})$ be an arc of type $(1,q-1)$. By $p_{(1,q-1),(1,q-1)}^{(1,q-2)}\neq0$ and Lemma \ref{jiben} (ii), there exists a vertex $z_{q-2}\in P_{(q-1,1),(1,q-2)}(z,z_{0})$. Since $(1,q-2)$ is pure, we assume that $(z_{0},z_{1},\ldots,z_{q-2})$ is a circuit consisting of arcs of the same type. Hence, $\wz{\partial}(z,z_{1})=(2,q-2)$. The fact that $\wz{\partial}(z,z_{0})=(1,q-1)$ and $\wz{\partial}(z_{0},z_{1})=(1,q-2)$ imply $p_{(1,q-1),(1,q-2)}^{(2,q-2)}\neq0$.

Suppose $p_{(1,q-1),(1,q-1)}^{(2,q-2)}\neq0$. Let $(y_{0},y_{1})$ and $(y_{1},y_{2})$ be arcs of type $(1,q-1)$ such that $\wz{\partial}(y_{0},y_{2})=(2,q-2)$. Since $p_{(1,q-1),(1,q-2)}^{(2,q-2)}\neq0$, there exists a vertex $y_{1}'\in P_{(1,q-1),(1,q-2)}(y_{0},y_{2})$. By Lemma \ref{jiben2} (i), one has $k_{1,q-1}=2$ and $p_{(1,q-1),(1,q-1)}^{(2,q-2)}=1$. Lemma \ref{jiben2} (ii) and (iii) imply that $p_{(1,q-1),(1,q-1)}^{(1,q-2)}=2$ and $(y_{1}',y_{1})$ is an arc of type $(1,q-1)$. Since $y_{0}\in P_{(q-1,1),(1,q-1)}(y_{1},y_{1}')$, from Lemma \ref{jiben2} (iv), we get $q=2$, a contradiction. Thus, (i) holds

(ii) Let $(x_{0},x_{1},\ldots,x_{q-1})$ be a circuit such that $\wz{\partial}(x_{q-1},x_{0})=(1,q-1)$. Suppose $\wz{\partial}(x_{0},x_{1})=(1,p-1)$ with $p\notin\{q,q-1\}$. It follows that $q>3$ and $\partial(x_{1},x_{q-1})=q-2$.

\textbf{Case 1}.~$\partial(x_{q-1},x_{1})=1$.

Since $x_{0}\in P_{(1,q-1),(1,p-1)}(x_{q-1},x_{1})$, there exists $y\in P_{(1,p-1),(1,q-1)}(x_{q-1},x_{1})$. By $p_{(1,q-1),(1,q-1)}^{(1,q-2)}\neq0$, we can pick a vertex $z\in P_{(1,q-1),(1,q-1)}(x_{q-1},x_{1})$. Note that $|\{x_{0},y,z\}|=3$. Since $(1,q-2)$ is pure, one gets $\wz{\partial}(x_{i},x_{i+1})=(1,q-2)$ for $1\leq i\leq q-2$, which implies $\{x_{0},y,z\}\subseteq\Gamma_{2,q-2}(x_{q-2})$, contrary to Lemma \ref{jiben2} (i).

\textbf{Case 2}.~$\partial(x_{q-1},x_{1})=2$.

Let $(z_{0}',z_{1}')$ and $(z_{1}',z_{2}')$ be arcs of type $(1,q-2)$ such that $\wz{\partial}(z_{0}',z_{2}')=(2,q-3)$. By $p_{(1,q-1),(1,q-1)}^{(1,q-2)}\neq0$, there exists a vertex $z'\in P_{(1,q-1),(1,q-1)}(z_{0}',z_{1}')$. Since $(1,q-2)$ is pure, one gets $\wz{\partial}(z',z_{2}')=(2,q-2)$. By $x_{0}\in P_{(1,q-1),(1,p-1)}(x_{q-1},x_{1})$, there exists a vertex $w\in P_{(1,q-1),(1,p-1)}(z',z_{2}')$, which implies $\wz{\partial}(z_{0}',w)=(2,q-2)$. Since $\wz{\partial}(z_{0}',z')=\wz{\partial}(z',w)=(1,q-1)$, we have $p_{(1,q-1),(1,q-1)}^{(2,q-2)}\neq0$, contrary to (i).

Note that $\wz{\partial}(x_{i},x_{i+1})=(1,q-1)$ or $(1,q-2)$ for $0\leq i\leq q-2$. If $\wz{\partial}(x_{i},x_{i+1})=(1,q-1)$ for each $i$, by $q\geq3$, then $\wz{\partial}(x_{0},x_{2})=(2,q-2)$, contrary to $p_{(1,q-1),(1,q-1)}^{(2,q-2)}=0$. Thus, (ii) holds.

(iii)~By Lemma \ref{jiben2} (ii), $k_{1,q-1}=2$. Since $k_{1,q-2}=1$, we get $p_{(1,q-1),(1,q-1)}^{(1,q-2)}=2$ from Lemma \ref{jiben} (v). Let $(w_{0}=w_{q-1},w_{1},\ldots,w_{q-2})$ be a circuit consisting of arcs of type $(1,q-2)$. Pick vertices $w'\in P_{(1,q-1),(1,q-1)}(w_{0},w_{1})$ and $w''\in P_{(1,q-1),(1,q-1)}(w_{1},w_{2})$ such that $\wz{\partial}(w',w'')\neq(1,q-2)$. Note that $(w'',w_{0})=q-2$. By (i), one has $q-2<\partial(w'',w')\leq1+\partial(w'',w_{0})=q-1$. Since $q\geq3$, we obtain $\partial(w',w'')=2$ from Lemma \ref{jiben2} (iv). The desired result follows.$\qed$

For any element $a$ in a residue class ring, we assume that $\hat{a}$ denotes the minimum nonnegative integer in $a$.

\begin{lemma}\label{claim 1}
Let $(1,h-1),(1,l-1)\in\wz{\partial}(\Gamma)$ and $v=\min\{j\mid p_{(1,h-1),(1,l-1)}^{(i,j)}\neq0\}$ with $h,l>2$. Suppose that $(1,l-1)$ is pure, or $p_{(1,l-1),(1,l-1)}^{(1,l-2)}\neq0$ and $(1,l-2)$ is pure. If $(\Gamma_{1,l-1})^{v}\cap\Gamma_{l-1,1}\Gamma_{h-1,1}\neq\emptyset$, then $h=l$ or $p_{(1,h-1),(1,h-1)}^{(1,l-2)}\neq0$.
\end{lemma}
\textbf{Proof.}~Let $(x_{0},x_{1},\ldots,x_{v+1})$ be a circuit of length $v+2$ such that $\wz{\partial}(x_{v+1},x_{0})=(1,h-1)$ and $\wz{\partial}(x_{i},x_{i+1})=(1,l-1)$ for $0\leq i\leq v$. Suppose that $h\neq l$.

\textbf{Case 1}.~$(1,l-1)$ is pure.

Note that $v+2>l$. Since $x_{0}\neq x_{l}$, by Lemma \ref{jiben2} (i), we have $k_{1,l-1}=2$. In view of $l>2$ and Lemma \ref{jiben3}, we get $\Delta_{l}\simeq{\rm Cay}(\mathbb{Z}_{2l},\{1,l+1\})$ or ${\rm Cay}(\mathbb{Z}_{l}\times\mathbb{Z}_{l},\{(1,0),(0,1)\})$.

\textbf{Case 1.1.} $\Delta_{l}\simeq{\rm Cay}(\mathbb{Z}_{2l},\{1,l+1\})$.

In view of Lemma \ref{jiben3} (i), we obtain $(\Gamma_{1,l-1})^{i}=\{\Gamma_{i,l-i}\}$ for $2\leq i\leq l-1$. Then $\wz{\partial}(x_{0},x_{l-1})=(l-1,1)$. If $v\geq l$, by Lemma \ref{jiben2} (iii), then $\wz{\partial}(x_{0},x_{l+1})=(1,l-1)$, contrary to the minimality of $v$; if $v=l-1$, by $x_{l-1}\in P_{(l-1,1),(1,l-1)}(x_{l},x_{0})$ and Lemma \ref{jiben2} (iv), then $h=2$, a contradiction.

\textbf{Case 1.2.} $\Delta_{l}\simeq{\rm Cay}(\mathbb{Z}_{l}\times\mathbb{Z}_{l},\{(1,0),(0,1)\})$.

Let $\tau$ be an isomorphism from ${\rm Cay}(\mathbb{Z}_{l}\times\mathbb{Z}_{l},\{(1,0),(0,1)\})$ to $\Delta_{l}$. Pick $\tau(a,b)\in\Gamma_{1,h-1}(\tau(0,0))$. Then $0\notin\{a,b\}$. Since $\tau(a,b)\in P_{(1,h-1),(l-\hat{a},\hat{a})}(\tau(0,0),\tau(0,b))$, we get $\tau(e+a,f+b)\in\Gamma_{1,h-1}(\tau(e,f))\cap\Gamma_{\hat{a},l-\hat{a}}(\tau(e,f+b))$ and $\tau(e+b,f+a)\in\Gamma_{1,h-1}(\tau(e,f))\cap\Gamma_{\hat{a},l-\hat{a}}(\tau(e+b,f))$ for each $(e,f)$. By $h\neq2$, one has $\hat{a}+\hat{b}\neq l$.

Suppose $a=-1$. Since $(\tau(0,0),\tau(1,0),\tau(1+a,b)=\tau(0,b),\tau(0,b+1),\ldots,\tau(0,l-1))$ is a circuit of length $l-\hat{b}+2$ containing arcs of types $(1,h-1)$ and $(1,l-1)$, we get $b=1$, contrary to $\hat{a}+\hat{b}\neq l$. Hence, $a\neq-1$. Similarly, $b\neq-1$. By $(\Gamma_{1,l-1})^{v}\cap\Gamma_{l-1,1}\Gamma_{h-1,1}\neq\emptyset$ and the minimality of $v$, one gets $v=\partial_{\Gamma}(\tau(a+1,b),\tau(0,0))=2l-\hat{a}-\hat{b}-1$. By $l-1\leq v$, we obtain $\hat{a}+\hat{b}<l$. Note that $(\tau(a+1,b),\tau(a+b+1,a+b),\tau(a+b+2,a+b),\ldots,\tau(0,a+b),\tau(0,a+b+1),\tau(0,a+b+2),\ldots,\tau(0,0))$ is a path. If $a+b=-1$, then $l+1-\hat{a}-\hat{b}\geq v$, contrary to $l>2$; if $a+b\neq -1$, then $2l-2\hat{a}-2\hat{b}\geq v$, contrary to $\hat{a}+\hat{b}>1$.

\textbf{Case 2}.~$p_{(1,l-1),(1,l-1)}^{(1,l-2)}\neq0$ and $(1,l-2)$ is pure.

Since $h\neq l$ and $h,l>2$, one has $v\geq2$. By the minimality of $v$, we obtain $\partial(x_{j},x_{j+2})=2$ for $0\leq j\leq v-1$. Lemma \ref{jiben2} (ii) implies $|(\Gamma_{1,l-1})^{2}|=2$ and $k_{1,l-1}=2$. If $|P_{(1,l-1),(1,l-1)}(x_{0},x_{2})|=2$, then there exists a vertex $x_{1}'\in P_{(1,l-1),(1,l-1)}(x_{0},x_{2})$ such that $\wz{\partial}(x_{1}',x_{3})=(1,l-2)$, contrary to the minimality of $v$. Then $|P_{(1,l-1),(1,l-1)}(x_{0},x_{2})|=1$. By Lemma \ref{jiben2} (iii),  $p_{(1,l-1),(1,l-1)}^{(1,l-2)}=2$. It follows from Lemma \ref{jiben} (i) and (v) that $k_{1,l-2}=1$. In view of Lemma \ref{mix-(2,q-2)'} (iii), we have $\wz{\partial}(x_{j},x_{j+2})=(2,l-1)$ for $0\leq j\leq v-1$. Hence, $v\geq l-1$. By Lemma \ref{mix-(2,q-2)'} (i) and Lemma \ref{jiben} (iii), we obtain $h\neq l-1$.

Let $(y_{0},y_{1},\ldots,y_{v+1})$ be a path consisting of arcs of type $(1,l-1)$ such that $\wz{\partial}(y_{j},y_{j+2})=(2,l-1)$ for $0\leq j\leq v-1$. Pick $x_{v+1}'$ and $y_{v+1}'$ such that $\Gamma_{1,l-1}(x_{v})=\{x_{v+1},x_{v+1}'\}$ and $\Gamma_{1,l-1}(y_{v})=\{y_{v+1},y_{v+1}'\}$. Then $\wz{\partial}(x_{v-1},x_{v+1}')=\wz{\partial}(y_{v-1},y_{v+1}')=(1,l-2)$. Since $k_{1,l-2}=1$, by Lemma \ref{jiben} (iii) and the inductive hypothesis, we have $\wz{\partial}(x_{0},x_{v+1}')=\wz{\partial}(y_{0},y_{v+1}')$, which implies $\wz{\partial}(x_{0},x_{v+1})=\wz{\partial}(y_{0},y_{v+1})$. Thus, $\wz{\partial}(x_{0},x_{v+1})$ only depends on $v$.

Since $(1,l-2)$ is pure and $k_{1,l-2}=1$, each $\Delta_{l-1}(x_{i})$ is a circuit of length $l-1$, denoted by $(x_{i}=x_{0,i},x_{1,i},\ldots,x_{l-2,i})$, where the first subscription of $x$ are taken modulo $l-1$. The fact that $p_{(1,l-1),(1,l-1)}^{(1,l-2)}=2$ implies that $\wz{\partial}(x_{a,b},x_{a,b+1})=\wz{\partial}(x_{a,b+1},x_{a+1,b})=(1,l-1)$ for any $a\in\{0,1,\ldots,l-2\}$ and $b\in\{0,1,\ldots,v\}$. By $k_{1,l-2}=1$, one gets $\wz{\partial}(x_{j,v-j+1},x_{j+2,v-j-1})=(2,l-1)$ for $0\leq j\leq v-1$. Since $\wz{\partial}(x_{0},x_{v+1})$ only depends on $v$, we obtain $\wz{\partial}(x_{0,v+1},x_{v+1,0})=\wz{\partial}(x_{0},x_{v+1})=(h-1,1)$. Let $r$ be the minimal nonnegative integer such that $r\equiv v+1$ (mod $l-1$). It suffices to show that $r=l-2$. Note that $(x_{0,0},x_{1,0},\ldots,x_{r,0}=x_{v+1,0},x_{0,v+1})$ is a circuit. By $h\neq2$, $r\neq0$. Since $h\neq l-1$ and $(1,l-2)$ is pure, one gets $r=l-2$.

This completes the proof of Lemma \ref{claim 1}.$\qed$

\begin{lemma}\label{claim 2}
Let $q\geq3$. If $(1,q-1)$ is pure, or $p_{(1,q-1),(1,q-1)}^{(1,q-2)}\neq0$ and $(1,q-2)$ is pure, then {\rm (i)} in Theorem \ref{Main1} is valid.
\end{lemma}
\textbf{Proof.}~Let $x_{0},x,x_{1}$ be vertices such that $\wz{\partial}(x_{0},x)=(1,s-1)$, $\wz{\partial}(x,x_{1})=(1,t-1)$ and $\wz{\partial}(x_{0},x_{1})=(1,q-1)$. By Lemma \ref{jiben2} (iv) and $s<t$, we have $s\neq q$. Suppose $t\neq q$. Observe that $p_{(1,q-1),(1,q-1)}^{(2,q-2)}\neq0$ or $p_{(1,q-1),(1,q-1)}^{(1,q-2)}\neq0$. Pick $x_{2}$ such that $\wz{\partial}(x_{1},x_{2})=(1,q-1)$ and $\partial(x_{2},x_{0})=q-2$.

We claim that $\wz{\partial}(x,x_{2})=(2,q-1)$. If $(1,q-1)$ is pure, by $q\notin\{s,t\}$, then our claim is valid. Suppose that $p_{(1,q-1),(1,q-1)}^{(1,q-2)}\neq0$ and $(1,q-2)$ is pure. It follows from Lemma \ref{mix-(2,q-2)'} (i) that $\wz{\partial}(x_{0},x_{2})=(1,q-2)$.

Suppose $s=q-1$. Since $x_{0}\in P_{(q-2,1),(1,q-2)}(x,x_{2})$ and $x\neq x_{2}$, by Lemma \ref{jiben2} (iv), we get $\wz{\partial}(x,x_{2})=(1,1)$ or $(2,2)$. In view of Lemma \ref{pure-g=2} and $t\neq q$, $\wz{\partial}(x,x_{2})=(2,2)$. Since $t\neq q-1$, from Lemma \ref{mix-(2,q-2)'} (ii), one has $\partial(x_{2},x)>q-2$. Hence, $q=3$, $s=2$ and $t=4$. By $x_{1}\in P_{(1,3),(1,2)}(x,x_{2})$, there exists a vertex $x_{1}'\in P_{(1,2),(1,3)}(x_{2},x)$, which implies $\partial(x_{1},x_{1}')=2$. It follows from Lemma \ref{mix-(2,q-2)'} (i) that $\partial(x_{1}',x_{1})=2$. By Lemma \ref{jiben2} (ii), we obtain $(\Gamma_{1,2})^{2}=\{\Gamma_{1,1},\Gamma_{2,2}\}$. Since $x\neq x_{2}$, from Lemma \ref{jiben2} (i), one has $k_{1,1}=2$. In view of Lemma \ref{jiben} (i) and (v), we get $p_{(1,2),(1,2)}^{(1,1)}=1$. By Lemma \ref{jiben2} (ii), $p_{(1,2),(1,2)}^{(2,2)}=2$. Hence, $\wz{\partial}(x,x_{1})=(1,2)$, a contradiction. Thus, $s\neq q-1$. Similarly, $t\neq q-1$.

Since $t\notin\{q-1,q\}$, by Lemma \ref{mix-(2,q-2)'} (ii), we get $q-1\leq\partial(x_{2},x)\leq1+\partial(x_{2},x_{0})=q-1$. The fact that $s\notin\{q-1,q\}$ and $\partial(x_{2},x_{0})=q-2$ imply $\partial(x,x_{2})=2$. Therefore, our claim is valid.

Since $x\in P_{(1,s-1),(1,t-1)}(x_{0},x_{1})$, there exists a vertex $x'\in P_{(1,t-1),(1,s-1)}(x_{0},x_{1})$. Similarly, $\wz{\partial}(x',x_{2})=(2,q-1)$. Since $x_{1}\in\Gamma_{1,t-1}(x)\cap\Gamma_{1,s-1}(x')\cap\Gamma_{q-1,1}(x_{2})$, there exist vertices $y_{1}'\in P_{(1,q-1),(1,t-1)}(x,x_{2})$ and $y_{1}''\in P_{(1,q-1),(1,s-1)}(x,x_{2})$. It follows from Lemma \ref{jiben2} (i) that $k_{1,q-1}=2$. Similarly, $y_{1}',y_{1}''\in\Gamma_{2,q-1}(x_{0})$. Then $\Gamma_{1,s-1}\Gamma_{1,q-1}=\{\Gamma_{2,q-1}\}$. By Lemma \ref{jiben} (i), we have $k_{1,s-1}=p_{(1,s-1),(1,q-1)}^{(2,q-1)}$. Since $x_{1}\in P_{(1,t-1),(1,q-1)}(x,x_{2})$, from Lemma \ref{jiben} (iv), one gets $k_{1,s-1}=1$. Similarly, $k_{1,t-1}=1$. In view of $p_{(1,s-1),(1,t-1)}^{(1,q-1)}\neq0$, we obtain $k_{1,q-1}=1$, a contradiction. Thus, $t=q$. By Lemma \ref{jiben2} (iv), one has $s=2$.$\qed$

\noindent\textbf{Proof of Theorem \ref{Main1}.} (ii) (a)$\Rightarrow$(b): By way of contradiction, we may assume that $q$ is the minimal integer such that $(1,q-1)$ is mixed and (b) does not hold. Since $(1,1)$ is pure, $q\geq3$. Pick a circuit $(x_{0},x_{1},\ldots,x_{q-1})$ such that $\wz{\partial}(x_{q-1},x_{0})=(1,q-1)$ and $\wz{\partial}(x_{0},x_{1})=(1,c-1)$ with $c<q$.

Suppose $\wz{\partial}(x_{i},x_{i+1})=(1,d-1)$ for some $i\in\{1,2,\ldots,q-2\}$ and $d\notin\{q,c\}$. Without loss of generality, we may assume $i=q-2$. Lemmas \ref{pure-g=2}, \ref{claim 2} and the minimality of $q$ imply $\wz{\partial}(x_{q-2},x_{0})=\wz{\partial}(x_{q-1},x_{1})=(2,q-2)$. Since $x_{0}\in P_{(1,q-1),(1,c-1)}(x_{q-1},x_{1})$, there exist vertices $z_{0}\in P_{(1,c-1),(1,q-1)}(x_{q-1},x_{1})$ and $z_{q-1}\in P_{(1,c-1),(1,q-1)}(x_{q-2},x_{0})$. In view of Lemma \ref{jiben2} (i), $k_{1,q-1}=2$. By Lemmas \ref{pure-g=2}, \ref{claim 2} and the minimality of $q$, we get $\wz{\partial}(z_{q-1},x_{1})=(2,q-2)$ and $\Gamma_{1,q-1}\Gamma_{1,c-1}=\{\Gamma_{2,q-2}\}$. It follows from Lemma \ref{jiben} (i) that $k_{1,c-1}=p_{(1,q-1),(1,c-1)}^{(2,q-2)}$. Since $x_{q-1}\in P_{(1,d-1),(1,q-1)}(x_{q-2},x_{0})$, by Lemma \ref{jiben} (iv), we obtain $k_{1,c-1}=1$. Similarly, $k_{1,d-1}=1$.

Since $k_{1,q-1}=2$, by Lemma \ref{jiben} (i) and Lemma \ref{jiben2} (i), one gets $\wz{\partial}(x_{j},x_{j+1})=(1,q'-1)$ for some $j\in\{1,2,\ldots,q-3\}$, and $k_{1,q'-1}=2$. Without loss of generality, we may assume $j=1$. It follows from Lemmas \ref{pure-g=2}, \ref{claim 2} and the minimality of $q$ that $\wz{\partial}(z_{0},x_{2})=(2,q-2)$. Since $x_{1}\in P_{(1,q-1),(1,q'-1)}(z_{0},x_{2})$, we have $x_{q-1}$ or $z_{q-1}\in P_{(1,q'-1),(1,q-1)}(x_{q-2},x_{0})$, a contradiction. Hence, $\wz{\partial}(x_{i},x_{i+1})=(1,q-1)$ or $(1,c-1)$ for each $i$.

Since $c<q$, by Lemmas \ref{pure-g=2} and \ref{claim 1}, we have $\wz{\partial}(x_{i},x_{i+1})=(1,q-1)$ for some $i\in\{1,2,\ldots,q-2\}$. Without loss of generality, we may assume $i=q-2$. Suppose $\partial(x_{q-2},x_{0})=2$. Then $\wz{\partial}(x_{q-2},x_{0})=\wz{\partial}(x_{q-1},x_{1})=(2,q-2)$. Since $x_{q-1}\in P_{(1,q-1),(1,q-1)}(x_{q-2},x_{0})$, there exists a vertex $x_{0}'\in P_{(1,q-1),(1,q-1)}(x_{q-1},x_{1})$, which implies $\wz{\partial}(x_{q-2},x_{0}')=(2,q-2)$ and $k_{1,q-1}=2$ from Lemma \ref{jiben2} (i). Hence, $(\Gamma_{1,q-1})^{2}=\{\Gamma_{2,q-2}\}$. By Lemma \ref{jiben2} (iii), we get $p_{(1,q-1),(1,q-1)}^{(2,q-2)}=2$ and $\wz{\partial}(x_{0},x_{1})=(1,q-1)$, a contradiction. Thus, $\wz{\partial}(x_{q-2},x_{0})=(1,q-2)$ and $p_{(1,q-1),(1,q-1)}^{(1,q-2)}\neq0$.

Note that $(1,q-2)$ is mixed. By the minimality of $q$, $p_{(1,q-2),(1,q-2)}^{(1,q-3)}\neq0$ and $(1,q-3)$ is pure. It follows from Lemma \ref{mix-(2,q-2)'} (ii) that the path $(x_{0},x_{1},\ldots,x_{q-2})$ contains an arc of type $(1,q-3)$. Hence, $c=q-2$ and $\wz{\partial}(x_{j},x_{j+1})=(1,q-3)$ for $0\leq j\leq q-3$. By Lemma \ref{pure-g=2}, we get $q>4$. In view of Lemma \ref{claim 1}, we obtain $p_{(1,q-2),(1,q-2)}^{(1,q-4)}\neq0$, a contradiction. Thus, our desired result holds.

(b)$\Rightarrow$(c): It is obvious.

(c)$\Rightarrow$(a): Suppose for the contrary that $(1,q-1)$ is pure. By Lemma \ref{jiben2} (ii), we have $|(\Gamma_{1,q-1})^{2}|=2$ and $k_{1,q-1}=2$. Lemma \ref{jiben3} implies that $p_{(1,q-1),(1,q-1)}^{(2,q-2)}=1$ and there exists an isomorphism $\tau$ from  ${\rm Cay}(\mathbb{Z}_{q}\times\mathbb{Z}_{q},\{(1,0),(0,1)\})$ to $\Delta_{q}$. It follows from Lemma \ref{jiben2} (iii) and Lemma \ref{jiben} (i), (v) that $k_{1,s-1}=1$. Observe that $(\tau(0,0),\tau(1,1),\ldots,\tau(-1,-1))$ is a circuit consisting of arcs of type $(1,s-1)$. Since $s\neq q$ from Lemma \ref{jiben2} (iv), $(1,s-1)$ is mixed. Then $p_{(1,s-1),(1,s-1)}^{(1,s-2)}\neq0$. By Lemma \ref{jiben} (i), we get $(\tau(1,1),\tau(3,3))\in\Gamma_{1,s-2}$ and $k_{1,s-2}=1$. Note that $(\tau(0,0),\tau(1,0),\tau(1,1),\tau(3,3),\tau(4,4),\ldots,\tau(-1,-1))$ is a circuit of length $q$ containing arcs of types $(1,q-1)$ and $(1,s-2)$, contrary to the fact that $(1,q-1)$ is pure. Thus, we have the assertion.

(i) follows by (ii) and Lemma \ref{claim 2}.

(iii) By Lemma \ref{jiben2} (iv), $s\neq q$. (ii) implies that $p_{(1,q-1),(1,q-1)}^{(1,q-2)}\neq0$ and $(1,q-2)$ is pure. In view of Lemma \ref{jiben2} (ii) and Lemma \ref{mix-(2,q-2)'} (iii), we only need to consider the case that $|(\Gamma_{1,q-1})^{2}|=2$ and $k_{1,q-2}=2$. Then $k_{1,q-1}=2$. By Lemma \ref{jiben} (i) and (v), we have $p_{(1,q-1),(1,q-1)}^{(1,q-2)}=1$. Suppose $s\neq q-1$. In view of Lemma \ref{jiben2} (iii), one gets $p_{(1,q-1),(1,q-1)}^{(1,s-1)}=2$ and $k_{1,s-1}=1$. Let $(x_{0},x),(x,x_{1})$ and $(x,x_{1}')$ be arcs of type $(1,q-1)$ such that $\wz{\partial}(x_{0},x_{1})=(1,s-1)$ and $\wz{\partial}(x_{0},x_{1}')=(1,q-2)$. Pick vertices $x_{2},z$ such that $\wz{\partial}(x_{1},x_{2})=\wz{\partial}(x,z)=(1,s-1)$. Since $p_{(1,q-1),(1,q-1)}^{(1,s-1)}=2$, we obtain $\wz{\partial}(x_{1},z)=\wz{\partial}(x_{1}',z)=\wz{\partial}(z,x_{2})=(1,q-1)$ and $\wz{\partial}(x_{1}',x_{2})=(1,q-2)$. The fact that $x_{1}'\in P_{(1,q-2),(1,q-2)}(x_{0},x_{2})$ and $k_{1,s-1}=1$ imply that $(1,s-1)$ is mixed. It follows from (ii) that $\wz{\partial}(x_{0},x_{2})=(1,s-2)$ and $(1,q-2)$ is mixed, contrary to the fact that $(1,q-2)$ is pure.$\qed$

\section{Configuration $C_{q,h}$}

In this section, we will discuss some useful properties of the configuration $C_{q,h}$.

\begin{lemma}\label{jiben4}
Suppose that $C_{q,h}$ exists. Then $k_{1,h-1}=1$, $k_{1,q-1}=2$, $(1,q-1)$ is pure and $\Delta_{q}\simeq{\rm Cay}(\mathbb{Z}_{2q},\{1,q+1\})$. Moreover, if $(1,q)$ is mixed, then $k_{1,q}=2$.
\end{lemma}
\textbf{Proof.}~Pick four distinct vertices $x,y,z,w$ such that $\wz{\partial}(x,w)=\wz{\partial}(y,w)=(1,q-1)$ and $\wz{\partial}(x,z)=\wz{\partial}(z,y)=(1,h-1)$. By Lemma \ref{jiben2} (i), $k_{1,q-1}=2$. In view of $h>2$ and Lemma \ref{jiben2} (iv), we have $|(\Gamma_{1,h-1})^{2}|=1$. Since $(\Gamma_{1,q-1})^{2}=\{\Gamma_{2,q-2}\}$, from Theorem \ref{Main1} (ii), $(1,q-1)$ is pure. Lemma \ref{jiben3} implies $\Delta_{q}\simeq{\rm Cay}(\mathbb{Z}_{2q},\{1,q+1\})$. So that there exists a vertex $w'\in P_{(1,q-1),(q-1,1)}(x,y)$ with $w'\neq w$. Write $\wz{\partial}(x,y)=\wz{f}$. By Lemma \ref{jiben} (i) and (v), one has $k_{\wz{f}}=1$. Since $(\Gamma_{1,h-1})^{2}=\{\Gamma_{\wz{f}}\}$, we get $k_{1,h-1}=1$. If $(1,q)$ is mixed, then $k_{1,q}=2$ from Theorem \ref{Main1} (ii) and Lemma \ref{jiben} (i).$\qed$

\begin{lemma}\label{h1}
Suppose that $C_{q,h}$ exists. The following hold:\vspace{-0.3cm}
\begin{itemize}
\item [{\rm(i)}] If $(1,h-1)$ is pure, then $h=4$.\vspace{-0.3cm}

\item [{\rm(ii)}] If $(1,h-1)$ is mixed, then $h=3$.
\end{itemize}
\end{lemma}
\textbf{Proof.}~Let $(x,z),(z,y)$ be two arcs of type $(1,h-1)$. Observe $P_{(1,q-1),(q-1,1)}(x,y)\neq\emptyset$. It follows from Lemma \ref{jiben2} (iv) that $\partial(x,y)=\partial(y,x)$. In view of Lemma \ref{jiben4}, one has $k_{1,h-1}=1$. If $(1,h-1)$ is pure, by Lemma \ref{jiben2} (ii), then $\wz{\partial}(x,y)=(2,2)$ and $h=4$; if $(1,h-1)$ is mixed, by Theorem \ref{Main1} (ii), then $\wz{\partial}(x,y)=(1,1)$ and $h=3$.$\qed$

\begin{lemma}\label{(1,h-1)_is_pure pre}
If $C_{q,h}$ exists, then $\Gamma_{1,q-1}\Gamma_{1,h-1}=\{\Gamma_{2,q}\}$ and $\Gamma_{q,2}\in\Gamma_{1,h-1}(\Gamma_{1,q-1})^{q-1}$.
\end{lemma}
\textbf{Proof.}~Pick four distinct vertices $x,y,z,w$ such that $\wz{\partial}(x,y)=\wz{\partial}(x,w)=(1,q-1)$ and $\wz{\partial}(y,z)=\wz{\partial}(z,w)=(1,h-1)$. By Lemma \ref{jiben4}, $(1,q-1)$ is pure and $k_{1,h-1}=1$. In view of Lemma \ref{jiben2} (ii), we have $|\Gamma_{1,q-1}\Gamma_{1,h-1}|=1$. It follows from Theorem \ref{Main1} (i) that $\partial(x,z)=2$. Note that $q-1\leq\partial(z,x)\leq1+\partial(w,x)=q$. It suffices to show that $\partial(z,x)=q$.

Assume the contrary, namely, there exists a path $(z=x_{0},x_{1},\ldots,x_{q-1}=x)$. Suppose that $\wz{\partial}(x_{i},x_{i+1})=(1,p-1)$ for some $i\in\{0,1,\ldots,q-2\}$ and $p\neq q$. Since $k_{1,h-1}=1$, we obtain $\wz{\partial}(z,x_{1})\neq(1,h-1)$. Hence, $p\neq h$. Without loss of generality, we may assume $i=q-2$. Since $(1,q-1)$ is pure, one has $\partial(y,x_{q-2})=q-1$, which implies $\partial(x_{q-2},y)=2$. By $x\in P_{(1,p-1),(1,q-1)}(x_{q-2},y)$ and Lemma \ref{jiben2} (i), we get $\wz{\partial}(w,z)=(1,p-1)$, contrary to $h\geq3$. Hence, $\wz{\partial}(x_{i},x_{i+1})=(1,q-1)$ for each $i$. It follows from Lemma \ref{jiben4} and Lemma \ref{jiben3} (i) that $\wz{\partial}(z,x)=(q-1,1)$, contrary to Lemma \ref{jiben2} (i).$\qed$

\begin{lemma}\label{(1,h-1)_is_mix pre}
If $(1,q-1)$ is mixed and $C_{q-1,h}$ exists, then $\Gamma_{1,q-1}\Gamma_{1,h-1}=\{\Gamma_{2,q}\}$ and $\Gamma_{q,2}\in\Gamma_{1,h-1}(\Gamma_{1,q-2})^{q-2}\Gamma_{1,q-1}$.
\end{lemma}
\textbf{Proof.}~Let $x,y,z$ be vertices such that $\wz{\partial}(x,y)=(1,q-1)$ and $\wz{\partial}(y,z)=(1,h-1)$. By Theorem \ref{Main1} (ii), we have $p_{(1,q-1),(1,q-1)}^{(1,q-2)}\neq0$. Pick a vertex $w\in P_{(q-1,1),(1,q-2)}(x,y)$. It follows from Theorem \ref{Main1} (i) that $\partial(x,z)=2$. Since $h\notin\{q,q-1\}$ from Lemma \ref{jiben4}, by Lemma \ref{mix-(2,q-2)'} (ii), one obtains $\partial(z,x)\geq q-1$. In view of Lemma \ref{(1,h-1)_is_pure pre}, we get $\wz{\partial}(w,z)=(2,q-1)$, which implies $\wz{\partial}(x,z)=(2,q)$ from Lemma \ref{jiben2} (iv). The desired results follow by Lemma \ref{(1,h-1)_is_pure pre}.$\qed$

\begin{lemma}\label{(1,g-1) (1,1)=(2,g)}
Suppose $q\geq3$ and $p_{(1,q-1),(1,1)}^{(2,s)}\neq0$. The following hold:\vspace{-0.3cm}
\begin{itemize}
\item [{\rm(i)}] If $(1,q-1)$ is pure, then $s=q$ and $\Gamma_{q,2}\in\Gamma_{1,1}(\Gamma_{1,q-1})^{q-1}$.\vspace{-0.3cm}

\item [{\rm(ii)}] If $(1,q-1)$ is mixed and $s=q$, then $\Gamma_{q,2}\in\Gamma_{1,1}\Gamma_{1,q-1}(\Gamma_{1,q-2})^{q-2}.$
\end{itemize}
\end{lemma}
\textbf{Proof.}~(i) Note that $s=q-1$ or $q$. Suppose for the contrary that $s=q-1$. Let $x_{q-1},x_{q},x_{0}$ be three vertices such that $\wz{\partial}(x_{q-1},x_{q})=(1,q-1)$, $\wz{\partial}(x_{q},x_{0})=(1,1)$ and $\wz{\partial}(x_{q-1},x_{0})=(2,q-1)$. Pick a path $(x_{0},x_{1},\ldots,x_{q-1})$.

\textbf{Case 1.} $\partial(x_{i+1},x_{i})\notin\{1,q-1\}$ for some $i\in\{0,1,\ldots,q-2\}$.

Without loss of generality, we may assume $\wz{\partial}(x_{q-2},x_{q-1})=(1,p-1)$ with $p\notin\{2,q\}$. Since $(1,q-1)$ is pure, we get $\wz{\partial}(x_{q-2},x_{q})=(2,q-1)$ from Theorem \ref{Main1} (i). In view of $x_{q}\in P_{(1,q-1),(1,1)}(x_{q-1},x_{0})$, there exists a vertex $x'\in P_{(1,1),(1,q-1)}(x_{q-2},x_{q})$, which implies $k_{1,q-1}=2$ by Lemma \ref{jiben2} (i). Since $(1,q-1)$ is pure, we have $\wz{\partial}(x',x_{0})=(2,q-1)$ and $\Gamma_{1,q-1}\Gamma_{1,1}=\{\Gamma_{2,q-1}\}$. It follows from Lemma \ref{jiben} (i) that $k_{1,1}=p_{(1,q-1),(1,1)}^{(2,q-1)}$. In view of $x_{q-1}\in P_{(1,p-1),(1,q-1)}(x_{q-2},x_{q})$ and Lemma \ref{jiben} (iv), we obtain $k_{1,1}=1$. By $x_{q}\in P_{(1,q-1),(q-1,1)}(x',x_{q-1})$ and Lemma \ref{jiben2} (iv), one gets $\partial(x',x_{q-1})=\partial(x_{q-1},x')$. Since $x_{q-2}\in P_{(1,1),(1,p-1)}(x',x_{q-1})$, we obtain $\wz{\partial}(x_{q-1},x_{q-2})=(1,p-1)$, contrary to $p\neq2$.

\textbf{Case 2.} $\partial(x_{i+1},x_{i})\in\{1,q-1\}$ for $0\leq i\leq q-2$.

Let $r-1$ be the number of arcs of type $(1,q-1)$ in the path $(x_{0},x_{1},\ldots,x_{q-1})$. Lemma \ref{pure-g=2} implies $r>1$. Without loss of generality, we may assume $\wz{\partial}(x_{j},x_{j+1})=(1,q-1)$ with $q-r\leq j\leq q-2$.  It follows from Theorem \ref{Main1} (ii) that $\wz{\partial}(x_{j},x_{j+2})=(2,q-2)$ or $(2,q-1)$ for each $j$.

Suppose $\wz{\partial}(x_{j},x_{j+2})=(2,q-1)$ for some $j$. It follows from Lemma \ref{jiben2} (ii) that $(\Gamma_{1,q-1})^{2}=\{\Gamma_{2,q-2},\Gamma_{2,q-1}\}$ and $k_{1,q-1}=2$. Lemma \ref{jiben3} implies $p_{(1,q-1),(1,q-1)}^{(2,q-2)}=1$. By Lemma \ref{jiben2} (iii), $p_{(1,q-1),(1,q-1)}^{(2,q-1)}=2$. Hence, $\wz{\partial}(x_{q},x_{0})=(1,q-1)$, a contradiction.

Suppose $\wz{\partial}(x_{j},x_{j+2})=(2,q-2)$ for each $j$. Since $\wz{\partial}(x_{q-1},x_{0})\neq(1,q-1)$, we have $r<q$ from Lemma \ref{jiben3}. Hence, $\wz{\partial}(x_{q-r},x_{q})=(r,q-r)$. By Lemma \ref{pure-g=2}, $r=\frac{q}{2}$. Since $\wz{\partial}(x_{0},x_{r})=(\frac{q}{2},\frac{q}{2})$, there exists a path $(y_{r}=x_{r},y_{r+1},\ldots,y_{q}=x_{0})$ consisting of arcs of type $(1,q-1)$. Then $(x_{0},x_{1},\ldots,x_{r}=y_{r},y_{r+1},\ldots,y_{q-1})$ is a circuit of length $q$ containing arcs of types $(1,q-1)$ and (1,1), a contradiction.

(ii) It is an immediate consequence of Theorem \ref{Main1} (ii).$\qed$

Let $A_{i,j}$ denote a matrix with rows and columns indexed by $V\Gamma$ such that $(A_{i,j})_{x,y}=1$ if $\wz{\partial}(x,y)=(i,j)$, and $(A_{i,j})_{x,y}=0$ otherwise.

\begin{lemma}\label{Gamma 2,q}
Suppose that $q>2$, $(1,q-1)$ is pure and $\Delta_{q}\simeq{\rm Cay}(\mathbb{Z}_{2q},\{1,q+1\})$. The following hold:\vspace{-0.3cm}
\begin{itemize}
\item [{\rm(i)}] If $(1,1)\in\wz{\partial}(\Gamma)$, then $A_{1,q-1}A_{1,1}=A_{1,q-1}$ or $A_{1,q-1}A_{1,1}=k_{1,1}A_{2,q}$.\vspace{-0.3cm}

\item [{\rm(ii)}] If $(1,q)$ is mixed, then $A_{1,q-1}A_{1,q}=2A_{2,q-1}$ and $(A_{1,q})^{2}=2A_{1,q-1}$.\vspace{-0.3cm}

\item [{\rm(iii)}] If $(1,q)$ is mixed and $A_{1,q-1}A_{1,1}=A_{1,q-1}$, then $A_{1,q}A_{1,1}=A_{1,q}$.\vspace{-0.3cm}

\item [{\rm(iv)}] If $(1,q)$ is mixed and $A_{1,q-1}A_{1,1}=k_{1,1}A_{2,q}$, then $A_{1,q}A_{1,1}=k_{1,1}A_{2,q+1}$.
\end{itemize}
\end{lemma}
\textbf{Proof.}~(i) Suppose $p_{(1,q-1),(1,1)}^{(1,q-1)}\neq0$. Since $\Delta_{q}\simeq{\rm Cay}(\mathbb{Z}_{2q},\{1,q+1\})$, we obtain $p_{(1,q-1),(q-1,1)}^{(1,1)}=2$. By Lemma \ref{jiben} (i) and (v), we get $k_{1,1}=1$, which implies $A_{1,q-1}A_{1,1}=A_{1,q-1}$. Suppose $p_{(1,q-1),(1,1)}^{(1,q-1)}=0$. By Theorem \ref{Main1} (i), Lemma \ref{pure-g=2} and Lemma \ref{(1,g-1) (1,1)=(2,g)} (i), we have $\Gamma_{1,q-1}\Gamma_{1,1}=\{\Gamma_{2,q}\}$, which implies $A_{1,q-1}A_{1,1}=k_{1,1}A_{2,q}$ from Lemma \ref{jiben} (i).

(ii)~By Theorem \ref{Main1} (ii), we get $p_{(1,q),(1,q)}^{(1,q-1)}\neq0$. Since $k_{1,q-1}=2$, from Lemma \ref{jiben} (i) and Lemma \ref{jiben2} (i), we have $k_{1,q}=2$.

Let $x,y,z,w$ be vertices such that $\wz{\partial}(x,y)=(1,q-1)$, $\wz{\partial}(y,z)=(1,q)$ and $w\in P_{(1,q-1),(q,1)}(y,z)$. By Lemma \ref{jiben3} (i), we have $\wz{\partial}(x,w)=(2,q-2)$. In view of Theorem \ref{Main1} (i), one gets $\wz{\partial}(x,z)=(2,q-1)$, which implies $A_{1,q-1}A_{1,q}=2A_{2,q-1}$ from Lemma \ref{jiben} (i) and Lemma \ref{jiben2} (i).

By Lemma \ref{jiben3} (i), there exists a vertex $y'\in P_{(1,q-1),(1,q-1)}(x,w)$ with $y\neq y'$. Since $p_{(1,q-1),(1,q)}^{(2,q-1)}=2$, one has $\wz{\partial}(y',z)=(1,q)$, which implies $(A_{1,q})^{2}=2A_{1,q-1}$ from Lemma \ref{jiben} (i).

(iii) By Lemma \ref{jiben} (i), we have $k_{1,1}=1$. Let $x_{0},x_{1},x_{2},x_{3}$ be vertices such that $\wz{\partial}(x_{0},x_{2})=(1,q-1)$, $x_{1}\in P_{(1,q-1),(1,1)}(x_{0},x_{2})$ and $x_{3}\in P_{(1,q),(1,q)}(x_{0},x_{2})$. It follows from (ii) that $\wz{\partial}(x_{3},x_{1})=(1,q)$. Since $x_{1}\in P_{(1,q),(1,1)}(x_{3},x_{2})$, by Lemma \ref{jiben} (i), we get $A_{1,q}A_{1,1}=A_{1,q}$.

(iv) Let $z_{0},z_{1},z_{2},z_{0}'$ be vertices such that $\wz{\partial}(z_{0},z_{1})=(1,q)$, $\wz{\partial}(z_{1},z_{2})=(1,1)$ and $z_{0}'\in P_{(q,1),(1,q-1)}(z_{0},z_{1})$. Since $A_{1,q-1}A_{1,1}=k_{1,1}A_{2,q}$, $\wz{\partial}(z_{0}',z_{2})=(2,q)$. In view of (ii), one has $\wz{\partial}(z_{0},z_{2})\neq(1,q)$, which implies $\partial(z_{0},z_{2})=2$ from Lemma \ref{pure-g=2} and Theorem \ref{Main1} (i). Since $\wz{\partial}(z_{0}',z_{2})=(2,q)$, by Lemma \ref{jiben2} (iv), we get $\partial(z_{2},z_{0})\neq q$. It follows from Lemma \ref{mix-(2,q-2)'} (ii) that $\partial(z_{2},z_{0})=q+1$. The desired result holds by Lemma \ref{jiben} (i) and Lemma \ref{jiben2} (i).$\qed$

\begin{lemma}\label{C_q,h holds k_1,1}
Suppose that $(1,1)\in\wz{\partial}(\Gamma)$ and $C_{q,h}$ exists. The following hold:\vspace{-0.3cm}
\begin{itemize}
\item [{\rm(i)}] $h=3$, $k_{1,1}=1$ and $A_{1,q-1}A_{1,1}=A_{1,q-1}$.\vspace{-0.3cm}

\item [{\rm(ii)}] If $(1,q)$ is mixed, then $A_{1,q}A_{1,1}=A_{1,q}$.
\end{itemize}
\end{lemma}
\textbf{Proof.}~(i) Let $x,y,z$ be vertices such that $\wz{\partial}(x,y)=(1,q-1)$ and $\wz{\partial}(y,z)=(1,1)$. Suppose $\wz{\partial}(x,z)=(2,q)$. Since $C_{q,h}$ exists, by Lemma \ref{(1,h-1)_is_pure pre}, there exists a vertex $w\in P_{(1,q-1),(1,h-1)}(x,z)$, which implies $\wz{\partial}(z,y)=(1,h-1)$, contrary to $h\geq3$. It follows from Lemma \ref{jiben4} and Lemma \ref{Gamma 2,q} (i) that $A_{1,q-1}A_{1,1}=A_{1,q-1}$ and $\wz{\partial}(x,z)=(1,q-1)$. Then $(\Gamma_{1,h-1})^{2}=\{\Gamma_{1,1}\}$ and $h=3$. By Lemma \ref{jiben} (i), $k_{1,1}=1$.

(ii) It is an immediate consequence of Lemma \ref{jiben4} and Lemma \ref{Gamma 2,q} (iii).$\qed$

\begin{prop}\label{h=h'}
If $C_{q,h}$ and $C_{q',h'}$ both exist, then $h=h'$.
\end{prop}
\textbf{Proof.}~If $(1,1)\in\wz{\partial}(\Gamma)$, by Lemma \ref{C_q,h holds k_1,1} (i), then $h=h'=3$; if $(1,1)\notin\wz{\partial}(\Gamma)$, by Theorem \ref{Main1} (ii) and Lemma \ref{h1}, then $h=h'=4$.$\qed$

\section{Proof of Proposition \ref{jiegou}}

We shall prove Proposition \ref{jiegou} by contradiction. Suppose that C1--C6 do not hold. Let $\mathscr{B}$ be the set consisting of $(p,p-1)$ and $(p-1,p)$ where $(1,p-1)$ is mixed, $\mathscr{C}=\{(p,q)\mid C_{p,q}~\textrm{or}~C_{q,p}~\textrm{exists}\}$ and $\mathscr{D}=\{(p,q)\mid (p,p-1)\in\mathscr{B}~\textrm{and}~C_{p-1,q}~\textrm{exists,}~\textrm{or}~(q,q-1)\in\mathscr{B}~\textrm{and}~C_{q-1,p}~\textrm{exists}\}$.

Suppose that $C_{q,h}$ exists for some $q$ and $h$. In view of Lemma \ref{jiben4}, $(1,q-1)$ is pure. If $(1,1)\in K$, from Lemma \ref{C_q,h holds k_1,1} (i), then $h=3$; if $(1,1)\notin K$, from Lemma \ref{h1} and Theorem \ref{Main1} (ii), then $h=4$. Since C1, C2 and C3 do not hold, by Proposition \ref{h=h'}, there exists $(1,p-1)\in K$ such that $p\neq2$ and $(q,p)\notin\mathscr{B}\cup\mathscr{C}\cup\mathscr{D}$. Suppose that $C_{t,h}$ does not exist for any $t$ and $h$. Since the valency of $\Gamma$ is more than $3$, we may assume that $(1,q-1)\in K$ with $q\neq2$. Since C4, C5 and C6 do not hold, from Theorem \ref{Main1} (ii), there exists $(1,p-1)\in K$ such that $p\neq2$ and $(q,p)\notin\mathscr{B}\cup\mathscr{C}\cup\mathscr{D}$.

We set
\begin{eqnarray}
l=\min\left\{r\mid p_{(1,i-1),(1,j-1)}^{(2,r)}\neq0,~i\neq j,~i,j\geq3,~(i,j)\notin\mathscr{B}\cup\mathscr{C}\cup\mathscr{D}\right\}.\nonumber
\end{eqnarray}
Without loss of generality, we may assume $p_{(1,q-1),(1,p-1)}^{(2,l)}\neq0$. By Lemma \ref{mix-(2,q-2)'} (ii) and Theorem \ref{Main1} (ii), one has $l\geq3$.

Choose vertices $x$, $y$ and $z$ with $\wz{\partial}(x,y)=(1,q-1)$, $\wz{\partial}(y,z)=(1,p-1)$ and $\wz{\partial}(x,z)=(2,l)$. Then there exists $y'$ such that $\wz{\partial}(x,y')=(1,p-1)$ and $\wz{\partial}(y',z)=(1,q-1)$.

The minimality of $l$ will be used many times in the sequel, so we will not refer to it every time for the sake of simplicity.
We will reach a contradiction under the following two separate cases:\vspace{-0.3cm}
\begin{itemize}
\item [{\rm A)}] There exists a shortest path from $z$ to $x$ containing an arc of type $(1,h-1)$ with $h\notin\{2,q,p\}$.\vspace{-0.3cm}

\item [{\rm B)}] Each arc of any shortest path from $z$ to $x$ is of type $(1,1)$, $(1,q-1)$ or $(1,p-1)$.
\end{itemize}

\subsection{The case A}

Without loss of generality, we may assume that $(z=x_{0},x_{1},\ldots,x_{l}=x)$ is a path such that $\wz{\partial}(x_{0},x_{1})=(1,h-1)$. For each $i$, write $h_{i}=\partial(x_{i+1},x_{i})+1$.

\begin{step}\label{C_t,h exists}
~{\rm Show that $C_{t,h}$ exists for some $t$.}\vspace{-0.3cm}
\end{step}

Assume the contrary, namely, $C_{t,h}$ does not exist for any $t$. Suppose that $(h,q),(h,p)\in\mathscr{B}\cup\mathscr{C}\cup\mathscr{D}$. Observe $\{(h,q),(h,p)\}\nsubseteq\mathscr{B}$. It follows from Proposition \ref{h=h'} that $\{(h,q),(h,p)\}\nsubseteq\mathscr{C}\cup\mathscr{D}$. Without loss of generality, we may assume $(h,q)\in\mathscr{B}$ and $(h,p)\in\mathscr{C}\cup\mathscr{D}$. If $C_{h,p}$ exists, by Lemma \ref{jiben4}, then $q=h+1$ and $(1,q-1)$ is mixed, contrary to $(q,p)\notin\mathscr{D}$. If $(1,h-1)$ is mixed and $C_{h-1,p}$ exists, then $C_{q,p}$ exists, a contradiction. Thus, $(h,q)$ or $(h,p)\notin\mathscr{B}\cup\mathscr{C}\cup\mathscr{D}$.

Without loss of generality, we may assume that $(h,p)\notin\mathscr{B}\cup\mathscr{C}\cup\mathscr{D}$. Theorem \ref{Main1} (i) implies $\wz{\partial}(y,x_{1})=(2,l)$. Since $z\in P_{(1,p-1),(1,h-1)}(y,x_{1})$, there exists a vertex $y''\in P_{(1,h-1),(1,p-1)}(x,z)$, which implies $k_{1,p-1}=2$ from Lemma \ref{jiben2} (i). By Theorem \ref{Main1} (i) again, we get $\partial(y',x_{1})=2$ and $\wz{\partial}(y'',x_{1})=(2,l)$. Then $\Gamma_{1,p-1}\Gamma_{1,h-1}=\{\Gamma_{2,l}\}$ and $k_{2,l}=2$. Since $p_{(1,q-1),(1,p-1)}^{(2,l)}\neq0$, from Lemma \ref{jiben} (i) and (iv), we obtain $k_{1,h-1}=p_{(1,p-1),(1,h-1)}^{(2,l)}$=1. By $k_{2,l}=2$, $\partial(x_{1},y')<l$. Hence, $(h,q)\in\mathscr{B}\cup\mathscr{C}\cup\mathscr{D}$. Since $k_{1,h-1}=1$, one gets $(h,q)\in\mathscr{B}$ from Lemma \ref{jiben4}. Suppose that $(1,q-1)$ is mixed. By Theorem \ref{Main1} (ii) and Lemma \ref{jiben} (ii), one has $p_{(1,q-1),(1,q-1)}^{(1,h-1)}=k_{1,q-1}$, which implies $\wz{\partial}(y,y'')=(1,q-1)$. Since $z\in P_{(1,p-1),(p-1,1)}(y'',y)$, from Lemma \ref{jiben2} (iv), one obtains $q=2$, a contradiction. Now suppose that $(1,h-1)$ is mixed. By Theorem \ref{Main1} (ii) again, $\wz{\partial}(y'',y)=(1,h-1)$. In view of Lemma \ref{jiben2} (iv), $h=2$, a contradiction. Thus, $C_{t,h}$ exists for some $t$.

\begin{step}\label{(q,h),(p,h)}
~{\rm Show that $\{(q,h),(p,h)\}\nsubseteq\mathscr{C}\cup\mathscr{D}$.}\vspace{-0.3cm}
\end{step}

Suppose for the contrary that $\{(q,h),(p,h)\}\subseteq\mathscr{C}\cup\mathscr{D}$. By Step \ref{C_t,h exists} and Lemma \ref{jiben4}, we have $k_{1,h-1}=1$. We conclude that $(x_{1},x_{2},\ldots,x_{l})$ consists of arcs of types $(1,q-1)$ and $(1,p-1)$.

By Proposition \ref{h=h'}, $C_{q,h}$ exists, or $C_{q-1,h}$ exists and $(1,q-1)$ is mixed. Suppose $h_{l-1}=2$. By Lemma \ref{C_q,h holds k_1,1} (i) or (ii), we have $\wz{\partial}(x_{l-1},y)=(1,q-1)$. Theorem \ref{Main1} (i) implies $\partial(x_{l-1},z)=2$, contrary to $\partial(z,x_{l-1})<l$. Then $h_{j}\neq2$ for $1\leq j\leq l-1$. Step \ref{C_t,h exists} and Proposition \ref{h=h'} imply that $h_{j}\in\{q,p,h\}$ for any $j$. Since $h\notin\{q,p\}$, one gets $l\geq4$ from Lemma \ref{mix-(2,q-2)'} (ii) and Theorem \ref{Main1} (ii). If $h_{j}=h$ for any $j$, by Lemma \ref{h1} and $k_{1,h-1}=1$, then $(1,h-1)$ is pure and $h=4$, which imply $z=x_{4}$, a contradiction. In the path $(x_{1},x_{2},\ldots,x_{l})$, without loss of generality, we may assume that the number of arcs of type $(1,p-1)$ is not less than the number of arcs of type $(1,q-1)$.

Without loss of generality, we may assume $h_{l-1}=p$. By Proposition \ref{h=h'} again, $C_{p,h}$ exists, or $C_{p-1,h}$ exists and $(1,p-1)$ is mixed. Suppose that $C_{p-1,h}$ exists and $(1,p-1)$ is mixed. Then $p\geq4$. In view of Lemma \ref{jiben4} and Proposition \ref{h=h'}, we obtain $(q,p-1)\notin\mathscr{B}\cup\mathscr{C}\cup\mathscr{D}$. Lemma \ref{Gamma 2,q} (ii) implies $\wz{\partial}(x_{l-1},y')=(1,p-2)$. It follows from Theorem \ref{Main1} (i) that $\partial(x_{l-1},z)=2$, contrary to $\partial(z,x_{l-1})<l$. Hence, $C_{p,h}$ exists. Suppose $h_{1}=h$. Then $\wz{\partial}(y,x_{2})=(1,p-1)$. By Theorem \ref{Main1} (i) again, one has $\partial(x,x_{2})=2$, contrary to $\partial(x_{2},x)<l$. Thus, $h_{j}\in\{q,p\}$ for $1\leq j\leq l-1$.

Suppose $h_{j}=p$ for each $j$. Since $C_{p,h}$ exists, by Lemma \ref{(1,h-1)_is_pure pre}, we get $\partial(x_{1},y)=p$, and so $l\geq p$. In view of Lemma \ref{jiben4} and Lemma \ref{jiben3} (i), one has $\wz{\partial}(x_{1},x_{p})=(p-1,1)$. Since $\wz{\partial}(y',z)=(1,q-1)$, we obtain $\wz{\partial}(x,x_{1})\neq(1,p-1)$ and $l>p$. Let $x_{l+1}=y'$. By Lemma \ref{jiben3} (i), one gets $\wz{\partial}(x_{1},x_{p+2})=(1,p-1)$. Then $x_{p+2}=y'$ and $(y',z,x_{1})$ is a circuit, a contradiction. Therefore, our conclusion is valid.

Without loss of generality, we may assume $h_{l-3}=q$ and $h_{l-2}=p$. Observe that $C_{p,h}$ exists and $(q,h)\in\mathscr{C}\cup\mathscr{D}$. From Lemma \ref{jiben4} and Proposition \ref{h=h'}, we get $k_{1,q-1}=k_{1,p-1}=2$ and there exists a vertex $x_{l-1}'\in P_{(1,p-1),(1,p-1)}(x_{l-2},x)$ with $x_{l-1}'\neq x_{l-1}$. Hence, $x_{l-1},x_{l-1}'\in\Gamma_{2,l}(x_{l-3})$. In view of Lemma \ref{jiben} (i) and Lemma \ref{jiben2} (i), we obtain $A_{1,q-1}A_{1,p-1}=2A_{2,l}$. Since $\wz{\partial}(x_{l-1},z)=\wz{\partial}(x,x_{1})=(3,l-1)$ and $x\in P_{(1,p-1),(2,l)}(x_{l-1},z)$, there exists a vertex $z''\in P_{(2,l),(1,p-1)}(x,x_{1})$, which implies $\wz{\partial}(y',z'')=(1,q-1)$. Then $\wz{\partial}(y',x_{1})=(2,l)$ and $\wz{\partial}(z,x_{1})=(1,p-1)$, contrary to $h\neq p$. The desired result follows.

\vspace{3ex}

In the following, we reach a contradiction based on the above discussion.

By Step \ref{(q,h),(p,h)}, we may assume $(p,h)\notin\mathscr{C}\cup\mathscr{D}$. It follows from Step \ref{C_t,h exists} and Lemma \ref{jiben4} that $k_{1,h-1}=1$. In view of Theorem \ref{Main1} (i), we have $\partial(y,x_{1})=\partial(y',x_{1})=2$.

\textbf{Case 1}. $\partial(x_{1},y)=l$.

Since $y'\in P_{(1,p-1),(1,q-1)}(x,z)$, there exists a vertex $z'\in P_{(1,p-1),(1,q-1)}(y,x_{1})$. It follows from Lemma \ref{jiben2} (i) that $k_{1,p-1}=2$. By Theorem \ref{Main1} (i), we have $\wz{\partial}(x,z')=(2,l)$, which implies $\Gamma_{1,q-1}\Gamma_{1,p-1}=\{\Gamma_{2,l}\}$. In view of Lemma \ref{jiben} (i), $k_{1,q-1}=p_{(1,q-1),(1,p-1)}^{(2,l)}$. Observe $z\in P_{(1,p-1),(1,h-1)}(y,x_{1})$. Lemma \ref{jiben} (iv) implies $k_{1,q-1}=1$. Since $k_{1,h-1}=1$, we obtain $\partial(x_{1},y')<l$. Note that $(q,h)\in\mathscr{B}\cup\mathscr{C}\cup\mathscr{D}$. From Lemma \ref{jiben4}, one has $(q,h)\in\mathscr{B}$. By Theorem \ref{Main1} (ii), we get $\wz{\partial}(z,z')=(1,q-1)$ or $\wz{\partial}(z',z)=(1,h-1)$, contrary to Lemma \ref{jiben2} (iv).

\textbf{Case 2}. $\partial(x_{1},y)<l$.

Note that $(p,h)\in\mathscr{B}$. By Lemma \ref{h1}, $(1,p-1)$ is mixed and $p=h+1=5$. Since $k_{1,3}=1$, one gets $k_{1,4}=p_{(1,4),(1,4)}^{(1,3)}$ from Theorem \ref{Main1} (ii) and Lemma \ref{jiben} (ii). If $\partial(x_{1},y')=l$, 
then there exists a vertex $w\in P_{(1,q-1),(1,4)}(y',x_{1})$, which implies $\wz{\partial}(z,w)=(1,4)$, contrary to Lemma \ref{jiben2} (iv). Hence, $\partial(x_{1},y')<l$.

Pick a vertex $w'\in \Gamma_{1,3}(y)$. Since $k_{1,4}=p_{(1,4),(1,4)}^{(1,3)}$, one has $\wz{\partial}(z,w')=\wz{\partial}(w',x_{1})=(1,4)$. The fact that $\partial(x_{1},y')<l$ implies $(q,4)\in\mathscr{C}\cup\mathscr{D}$. By Proposition \ref{h=h'}, $C_{q,4}$ exists, or $C_{q-1,4}$ exists and $(1,q-1)$ is mixed.  In view of Lemma \ref{(1,h-1)_is_pure pre} or \ref{(1,h-1)_is_mix pre}, we get $q=\partial(x_{1},y')<\partial(z,x)\leq1+\partial(w',x)=q+1$. Thus, $l=q+1$.

Suppose that $C_{q,4}$ exists. Pick a vertex $x_{2}'\in P_{(1,3),(4,1)}(w',x_{1})$. Then $\wz{\partial}(x,x_{2}')=(1,q-1)$. By Lemma \ref{jiben4}, there exists a circuit $(x_{2}',x_{3}',\ldots,x_{l}'=x)$ consisting of arcs of type $(1,q-1)$. Since $(z,x_{1},x_{2}',x_{3}',\ldots,x_{l}')$ is a shortest path, $\wz{\partial}(x_{1},x_{3}')=(2,l)$. It follows that $\wz{\partial}(x_{l-1}',z)=\wz{\partial}(x,x_{1})=(3,l-1)$. The fact that $x\in P_{(1,q-1),(2,l)}(x_{l-1}',z)$ and $\partial(x_{1},y)<l$ imply $\wz{\partial}(x_{2}',x_{1})=(2,l)$, a contradiction.

Suppose that $C_{q-1,4}$ exists and $(1,q-1)$ is mixed. Since $(1,4)$ is mixed, by Lemma \ref{jiben4}, we obtain $q\geq7$. It follows from Lemma \ref{(1,h-1)_is_mix pre} that there exists a vertex $y_{1}$ such that $\wz{\partial}(w',y_{1})=(1,q-2)$ and $\partial(y_{1},x)=q-1$. By Proposition \ref{h=h'}, we have $(5,q-1)\notin\mathscr{B}\cup\mathscr{C}\cup\mathscr{D}$. Theorem \ref{Main1} (i) implies $\wz{\partial}(z,y_{1})=(2,l)=(2,q+1)$. In view of $w'\in P_{(1,4),(1,q-2)}(z,y_{1})$, there exists a vertex $y''\in P_{(1,4),(1,q-2)}(x,z)$. By Lemma \ref{jiben4} and Lemma \ref{Gamma 2,q} (ii), we get $\wz{\partial}(y'',y')=(1,q-1)$, contrary to Lemma \ref{jiben2} (iv).

By the above discussion, we finish the proof of Proposition \ref{jiegou} for the case A.

\subsection{The case B}

Let $(z=x_{0},x_{1},\ldots,x_{l}=x)$ be a path. For each $i$, write $h_{i}=\partial(x_{i+1},x_{i})+1$. Note that $h_{i}\in\{2,q,p\}$.

\begin{step}\label{(1,g-1) or (1,q-1)}
~{\rm Show that $|\{i\mid h_{i}\neq2,~0\leq i\leq l-1\}|\geq2$ and $|\{h_{i}\mid0\leq i\leq l-1\}|\geq2$.}\vspace{-0.3cm}
\end{step}

Suppose that $h_{j}=2$ for $0\leq j\leq l-2$. It follows from Lemma \ref{pure-g=2} that $l=3$ or $4$. In view of Lemma \ref{jiben2} (i), we obtain $k_{1,1}=2$. By Lemma \ref{pure-g=2}, $\wz{\partial}(x_{0},x_{2})=(2,2)$ and $\Delta_{2}$ is not isomorphic to $C_{3}$. If $(1,3)$ is pure, then there exists a vertex $x_{1}'\in P_{(1,3),(1,3)}(x_{2},x_{0})$, a contradiction. Then $(1,3)\notin\wz{\partial}(\Gamma)$ or $(1,3)$ is mixed. Since $(q,p)\notin\mathscr{B}$, we get $\{q,p\}\neq\{3,4\}$ from Theorem \ref{Main1} (ii). By Lemma \ref{mix-(2,q-2)'} (ii), one has $l=4$ and $5\in\{q,p\}$. Lemma \ref{pure-g=2} and Theorem \ref{Main1} (i) imply $\partial(y,x_{1})=\partial(y',x_{1})=2$. Since $\partial(x_{1},y)\leq4$ and $\partial(x_{1},y')\leq4$, $(1,4)$ is mixed from Lemma \ref{(1,g-1) (1,1)=(2,g)} (i). By Theorem \ref{Main1} (ii), $(1,3)$ is pure, a contradiction. Therefore, the first statement is valid. The second statement follows from Lemma \ref{claim 1} and Theorem \ref{Main1} (iii).

\begin{step}\label{step}
~{\rm Show that $k_{1,q-2}=1$ if $|\{i\mid h_{i}=q\}|\geq2$ and $(1,q-1)$ is mixed.}\vspace{-0.3cm}
\end{step}

Without loss of generality, we may assume that $h_{l-2}=h_{l-1}=q$. Note that $\partial(x_{l-2},x)=2$. By Theorem \ref{Main1} (ii) and Lemma \ref{jiben2} (ii), we have $|(\Gamma_{1,q-1})^{2}|=2$ and $k_{1,q-1}=2$. Suppose $p_{(1,q-1),(1,q-1)}^{(1,q-2)}=1$. It follows from Lemma \ref{jiben2} (iii) that there exists a vertex $x_{l-1}'\in P_{(1,q-1),(1,q-1)}(x_{l-2},x)$ such that $\wz{\partial}(x_{l-1}',y)=(1,q-2)$. Pick a vertex $x'\in P_{(1,q-2),(1,q-1)}(x_{l-2},y)$. By Theorem \ref{Main1} (i), one gets $\partial(x',z)=2$, contrary to $\partial(z,x')<l$. Hence, $p_{(1,q-1),(1,q-1)}^{(1,q-2)}=2$. In view of Lemma \ref{jiben} (i) and (v), we obtain $k_{1,q-2}=1$.

\begin{step}\label{partial neq 1}
~{\rm Show that $\partial(x_{l-1},z)\geq2$.}\vspace{-0.3cm}
\end{step}

Suppose for the contrary that $\partial(x_{l-1},z)=1$.

\textbf{Case 1.} $(1,l-1)$ is mixed.

By Theorem \ref{Main1} (ii), $(1,l-2)$ is pure and there exists a vertex $x_{l-2}'$ such that $\wz{\partial}(x_{l-2}',x_{l-1})=(1,l-1)$ and $\wz{\partial}(x_{l-2}',z)=(1,l-2)$. Observe that $l-1,l\in\{2,q,p\}$. Since $(q,p)\notin\mathscr{B}$, $l=3$. From Lemma \ref{mix-(2,q-2)'} (ii), $\{q,p\}=\{3,4\}$. Without loss of generality, we may assume $p=4$. By Theorem \ref{Main1} (ii), $(1,3)$ is pure. In view of $\partial(x_{l-2}',y)\leq3$ and Lemma \ref{(1,g-1) (1,1)=(2,g)} (i), we get $\partial(y,x_{l-2}')=1$. It follows from Lemma \ref{pure-g=2} and Theorem \ref{Main1} (i) that $\wz{\partial}(y,x_{l-2}')=(1,3)$ and $\partial(x,x_{l-2}')=2$, contrary to $\partial(x_{l-2}',x)<l$.

\textbf{Case 2.} $(1,l-1)$ is pure.

Observe that $l\in\{q,p\}$.

\textbf{Case 2.1.} $h_{l-1}\neq2$.

Without loss of generality, we may assume $h_{l-1}=q$.  By $l\geq3$ and Lemma \ref{jiben2} (iv), one has $l=p$, which implies $h_{j}=p$ for $0\leq j\leq l-2$. In view of Lemma \ref{jiben2} (i), $k_{1,p-1}=2$.

We claim that $k_{1,q-1}=2$ and there exists $z'\in P_{(1,p-1),(1,p-1)}(y,x_{1})\setminus\{z\}$. Lemma \ref{jiben3} implies $\Delta_{p}\simeq{\rm Cay}(\mathbb{Z}_{2p},\{1,p+1\})$ or ${\rm Cay}(\mathbb{Z}_{p}\times\mathbb{Z}_{p},\{(1,0),(0,1)\})$. Suppose $\Delta_{p}\simeq{\rm Cay}(\mathbb{Z}_{2p},\{1,p+1\})$. Since $C_{p,q}$ does not exist, from Lemma \ref{jiben2} (ii), we have $|(\Gamma_{1,q-1})^{2}|=2$ and $k_{1,q-1}=2$. It follows from Lemma \ref{jiben3} (i) that the claim is valid. Suppose $\Delta_{p}\simeq{\rm Cay}(\mathbb{Z}_{p}\times\mathbb{Z}_{p},\{(1,0),(0,1)\})$. Then $|P_{(1,p-1),(p-1,1)}(x_{l-1},y)|=1$. Lemma \ref{jiben} (v) and Lemma \ref{jiben2} (i) imply $k_{\wz{\partial}(x_{l-1},y)}=2$. By Lemma \ref{jiben} (i), $k_{1,q-1}=2$. Since $\wz{\partial}(x_{l-1},x_{1})=(2,p-2)$, one has $\wz{\partial}(y,x_{1})\neq(2,p-2)$ from Lemma \ref{jiben3} (ii). In view of Lemma \ref{jiben2} (iii), the claim is valid.

By Theorem \ref{Main1} (i), $\wz{\partial}(x,z')=(2,l)$. In view of Lemma \ref{jiben} (i) and Lemma \ref{jiben2} (i), one gets $A_{1,q-1}A_{1,p-1}=2A_{2,l}$, which implies that $\wz{\partial}(x,x_{l-1})=(1,q-1)$, contrary to $q\neq2$.

\textbf{Case 2.2.} $h_{l-1}=2$.

Without loss of generality, we may assume $l=q$. By $x_{l-1}\neq y'$ and Lemma \ref{jiben2} (i), one gets $k_{1,q-1}=2$. Since $z\in P_{(1,q-1),(q-1,1)}(x_{l-1},y')$, we have $\wz{\partial}(y',x_{l-1})=(2,2)$ from Lemma \ref{jiben2} (iv) and Lemma \ref{pure-g=2}. In view of $x\in P_{(1,1),(1,p-1)}(x_{l-1},y')$, there exists a vertex $x''\in P_{(1,p-1),(1,1)}(y',x_{l-1})$. Then $x''\neq x$ and $k_{1,1}=2$. The fact that $\wz{\partial}(x_{l-1},y')=(2,2)$ implies $p_{(1,q-1),(q-1,1)}^{(1,1)}=p_{(1,1),(1,q-1)}^{(1,q-1)}=0$. Since $(1,q-1)$ is pure, by Lemma \ref{pure-g=2}, Theorem \ref{Main1} (i) and Lemma \ref{(1,g-1) (1,1)=(2,g)} (i), we obtain $A_{1,1}A_{1,q-1}=2A_{2,q}$. Hence, $y'\in P_{(1,1),(1,q-1)}(x'',z)$, a contradiction.

\begin{step}\label{pre-main3 1}
~{\rm Show that $p_{(1,s-1),(1,1)}^{(2,l-1)}\neq0$ for some $s>2$ if $\partial(x_{l-1},z)=2$.}\vspace{-0.3cm}
\end{step}

Pick a path $(x_{l-1},w,z)$ such that $\wz{\partial}(x_{l-1},w)=(1,s-1)$, $\wz{\partial}(w,z)=(1,t-1)$ and $s\geq t$. By Step \ref{(1,g-1) or (1,q-1)}, we may assume $h_{0}=q$. If $t=2$, then $s>2$ since $l\geq3$, and the desired result holds. Suppose $t\neq2$. Since $(q,p)\notin\mathscr{C}\cup\mathscr{D}$, from Lemmas \ref{(1,h-1)_is_pure pre} and \ref{(1,h-1)_is_mix pre}, we have $t=s$ or $(s,t)\in\mathscr{B}$.

\textbf{Case 1.} $s=t$.

Since $\partial(x_{1},w)<l$, by Theorem \ref{Main1} (i), one has $t=q$ or $(t,q)\in\mathscr{B}\cup\mathscr{C}\cup\mathscr{D}$.

\textbf{Case 1.1.} $t=q$.

Suppose $\wz{\partial}(x_{l-1},y')=(1,q-1)$. By Theorem \ref{Main1} (i), $h_{l-1}=p$. It follows from Theorem \ref{Main1} (ii) and (iii) that $(1,p-1)$ is mixed and $q=p-1$, contrary to $(q,p)\notin\mathscr{B}$. Since $w\neq y'$, by Lemma \ref{jiben2} (i), one gets $k_{1,q-1}=2$ and $p_{(1,q-1),(1,q-1)}^{(2,l-1)}=1$.

\textbf{Case 1.1.1.} $(1,q-1)$ is pure.

By Lemma \ref{jiben3}, one has $|(\Gamma_{1,q-1})^{2}|=2$ and $p_{(1,q-1),(1,q-1)}^{(2,q-2)}=1$. In view of Lemma \ref{jiben2} (iii), we get $\wz{\partial}(x_{l-1},z)=(2,q-2)$, which implies $l=q-1$ and $\wz{\partial}(w,x_{1})=(2,q-2)$. It follows from Theorem \ref{Main1} (ii) that $\wz{\partial}(y',x_{1})=(2,q-1)$. By Lemma \ref{jiben2} (iii) again, we obtain $p_{(1,q-1),(1,q-1)}^{(2,q-1)}=2$. Then $\wz{\partial}(x,y')=(1,q-1)$, a contradiction.

\textbf{Case 1.1.2.} $(1,q-1)$ is mixed.

By Theorem \ref{Main1} (ii), $(1,q-2)$ is pure. Since $p_{(1,q-1),(1,q-1)}^{(2,l-1)}=1$, one gets $p_{(1,q-1),(1,q-1)}^{(1,q-2)}=2$ from Lemma \ref{jiben2} (ii) and (iii). In view of Lemma \ref{jiben} (i) and (v), one has $k_{1,q-2}=1$. Lemma \ref{mix-(2,q-2)'} (iii) implies $\wz{\partial}(x_{l-1},z)=(2,q-1)$ and $l=q$. Pick vertices $x_{l-2}'\in P_{(q-1,1),(1,q-2)}(x_{l-1},w)$ and $x_{1}'\in P_{(1,q-2),(q-1,1)}(w,z)$. Note that $\wz{\partial}(x_{l-2}',x_{1}')\neq(2,q-3)$. By $k_{1,q-2}=1$, we obtain $l=q=3$. In view of Lemma \ref{mix-(2,q-2)'} (ii) and Theorem \ref{Main1} (ii), we get $p=4$, which implies that $(1,3)$ is pure. Since $\wz{\partial}(x_{l-1},z)=(2,2)$, there exists a vertex $x'\in P_{(1,3),(1,3)}(z,x_{l-1})$. Then  $(x_{l-1},w,z,x')$ is a circuit containing arcs of types $(1,2)$ and $(1,3)$, a contradiction.

\textbf{Case 1.2.} $(t,q)\in\mathscr{B}$ and $t=q-1$.

Note that $q>3$. Theorem \ref{Main1} (ii) implies that $(1,q-2)$ is pure. Since $q-1\notin\{2,q,p\}$, we have $\wz{\partial}(x_{l-1},z)\neq(2,q-3)$. In view of Lemma \ref{jiben2} (ii), we get $|(\Gamma_{1,q-2})^{2}|=2$ and $k_{1,q-2}=2$, which imply $p_{(1,q-2),(1,q-2)}^{(2,q-3)}=1$ from Lemma \ref{jiben3}. By Lemma \ref{jiben2} (iii), we obtain $p_{(1,q-2),(1,q-2)}^{(2,l-1)}=2$. In view of Lemma \ref{jiben} (i) and (v), one gets $k_{2,l-1}=1$.

Since $k_{1,q-2}=2$ and $(q,p)\notin\mathscr{D}$, by Lemma \ref{jiben4}, we have $(q-1,p)\notin\mathscr{B}\cup\mathscr{C}\cup\mathscr{D}$. It follows from $\partial(w,x_{l-2})<l$ and Theorem \ref{Main1} (i) that $h_{l-2}\neq p$. Hence, $h_{i}\neq p$ for $0\leq i\leq l-2$. In view of Step \ref{step}, we get $|\{j\mid h_{j}=q\}|<2$. By Step \ref{(1,g-1) or (1,q-1)}, one has $h_{l-1}=p$.

Since $\partial(x_{1},x)<l$ and $(q-1,p)\notin\mathscr{B}\cup\mathscr{C}\cup\mathscr{D}$, we have $\wz{\partial}(y',x_{1})\neq(1,q-2)$ from Theorem \ref{Main1} (i). In view of Theorem \ref{Main1} (ii), we get $p_{(1,q-1),(1,q-1)}^{(1,q-2)}\neq0$, which implies $|(\Gamma_{1,q-1})^{2}|=2$ and $k_{1,q-1}=2$ from Lemma \ref{jiben2} (ii). Since $k_{1,q-2}=2$,  by Lemma \ref{jiben} (i) and (v), we get $p_{(1,q-1),(1,q-1)}^{(1,q-2)}=1$.  It follows from Lemma \ref{jiben2} (iii) that there exists a vertex $z'\in P_{(1,q-1),(1,q-1)}(y',x_{1})\setminus\{z\}$. In view of Theorem \ref{Main1} (i), we obtain $z'\in\Gamma_{2,l}(x)$ and $\Gamma_{1,p-1}\Gamma_{1,q-1}=\{\Gamma_{2,l}\}$. Since $x\in P_{(1,p-1),(2,l)}(x_{l-1},z)$ and $k_{2,l-1}=1$, we obtain $k_{1,p-1}=2$ from Lemma \ref{jiben} (i) and Lemma \ref{jiben2} (i). Hence, $p_{(1,p-1),(1,q-1)}^{(2,l)}=2$ and there exists a vertex $y''\in P_{(1,p-1),(1,q-1)}(x,z)$ such that $\wz{\partial}(y'',x_{1})=(1,q-2)$. By Theorem \ref{Main1} (i), one has $\partial(x,x_{1})=2$, contrary to $\partial(x_{1},x)<l$.

\textbf{Case 1.3.} $(t,q)\in\mathscr{B}$ and $t=q+1$.

Since $(1,q)$ is mixed, $(1,q-1)$ is pure and $p_{(1,q),(1,q)}^{(1,q-1)}\neq0$ from Theorem \ref{Main1} (ii). By $\partial(x_{l-1},z)=2$ and Lemma \ref{jiben2} (ii), we have $|(\Gamma_{1,q})^{2}|=2$ and $k_{1,q}=2$. If $p_{(1,q),(1,q)}^{(2,l-1)}=2$, then there exists a vertex $w'\in P_{(1,q),(1,q)}(x_{l-1},z)$ such that $\wz{\partial}(y',w')=(1,q)$; if $p_{(1,q),(1,q)}^{(2,l-1)}=1$, by Lemma \ref{jiben2} (iii), then $p_{(1,q),(1,q)}^{(1,q-1)}=2$ and $\wz{\partial}(y',w)=(1,q)$. Without loss of generality, we may assume $\wz{\partial}(y',w)=(1,q)$.

By Theorem \ref{Main1} (i), we have $\partial(x,w)=2$. Since $k_{1,q}=2$ and $C_{q,p}$ does not exist, we obtain $(q+1,p)\notin\mathscr{B}\cup\mathscr{C}\cup\mathscr{D}$ from Lemma \ref{jiben4}. It follows that $l\leq\partial(w,x)\leq\partial(w,x_{l-1})+1=q+1$. In view of Lemma \ref{mix-(2,q-2)'} (i), one gets $l-1=\partial(z,x_{l-1})>q-1$. Then $\partial(w,x)=l$. Since $x\in P_{(l,2),(2,l)}(w,z)$, by Lemma \ref{jiben2} (iv), one has $q=1$, a contradiction.

\textbf{Case 1.4.} $(t,q)\in\mathscr{C}\cup\mathscr{D}$.

Suppose that $(1,t-1)$ is pure. By Lemma \ref{jiben4}, $k_{1,t-1}=1$ or $\Delta_{t}\simeq{\rm Cay}(\mathbb{Z}_{2t},\{1,t+1\})$, which implies $\wz{\partial}(x_{l-1},z)=(2,t-2)$. Hence, $t=p$, a contradiction.

Suppose that $(1,t-1)$ is mixed. It follows from Lemma \ref{jiben4} that $k_{1,t-1}=1$ or $C_{t-1,q}$ exists. If $k_{1,t-1}=1$, by Theorem \ref{Main1} (ii), then $\partial(x_{l-1},z)=1$, a contradiction; if $C_{t-1,q}$ exists, by Lemma \ref{jiben4} and Lemma \ref{Gamma 2,q} (ii), then $\partial(x_{l-1},z)=1$, a contradiction.

\textbf{Case 2.} $(s,t)\in\mathscr{B}$.

Note that $(1,s-1)$ is mixed and $s=t+1$. By Theorem \ref{Main1} (ii), $(1,t-1)$ is pure. Since $\partial(x_{1},w)<l$, from Theorem \ref{Main1} (i), one has $t=q$ or $(t,q)\in\mathscr{B}\cup\mathscr{C}\cup\mathscr{D}$.

\textbf{Case 2.1.} $t=q$.

Note that $s=q+1$ and $l-1=\partial(z,x_{l-1})\geq q-1$. Since $(1,q)$ is mixed, by Theorem \ref{Main1} (ii), one has $p_{(1,q),(1,q)}^{(1,q-1)}\neq0$. Pick a vertex $x_{l-2}'\in P_{(q,1),(1,q-1)}(x_{l-1},w)$. The fact that $q+1\notin\{2,q,p\}$ implies $\wz{\partial}(x_{l-2}',z)\neq(2,q-2)$. By Lemma \ref{jiben2} (ii), we get $|(\Gamma_{1,q-1})^{2}|=2$ and $k_{1,q-1}=2$. In view of Lemma \ref{jiben3}, we get $p_{(1,q-1),(1,q-1)}^{(2,q-2)}=1$. Since $\wz{\partial}(x_{l-2}',z)\neq(2,q-2)$, we obtain $\wz{\partial}(x_{l-2}',y')=(1,q-1)$ from Lemma \ref{jiben2} (iii).

Note that $p_{(1,q),(1,q)}^{(1,q-1)}\neq0$ and $k_{1,q-1}=2$. By Lemma \ref{jiben} (i) and Lemma \ref{jiben2} (i), we obtain $k_{1,q}=2$. Since $C_{q,p}$ does not exist, $(q+1,p)\notin\mathscr{B}\cup\mathscr{C}\cup\mathscr{D}$ from Lemma \ref{jiben4}. Since $x\in P_{(1,h_{l-1}-1),(1,p-1)}(x_{l-1},y')$, there exists a vertex $x'\in P_{(1,p-1),(1,h_{l-1}-1)}(x_{l-1},y')$. In view of Theorem \ref{Main1} (i), we have $\partial(x_{l-2}',x')=2$. Since $q\leq l$ and $l\leq\partial(x',x_{l-2}')\leq1+\partial(y',x_{l-2}')=q$, we get $\wz{\partial}(x_{l-2}',x')=(2,l)$.

Suppose $w=y'$. Since $x\in P_{(1,h_{l-1}-1),(1,p-1)}(x_{l-1},w)$ and $p\neq q+1$, by Theorem \ref{Main1} (i), we have $h_{l-1}=p$. Theorem \ref{Main1} (ii) and (iii) imply that $(1,q)$ is pure, a contradiction. Suppose $w\neq y'$. By $p_{(1,q-1),(1,p-1)}^{(2,l)}\neq0$ and $p>2$, we get $w\in P_{(1,q-1),(1,p-1)}(x_{l-2}',x')$. Since $x'\in P_{(1,p-1),(p-1,1)}(x_{l-1},w)$, from Lemma \ref{jiben2} (iv), we obtain $q=1$, a contradiction.

\textbf{Case 2.2.} $(t,q)\in\mathscr{B}$.

Note that $t=q-1$ and $s=q$. By Theorem \ref{Main1} (ii), $p_{(1,q-1),(1,q-1)}^{(1,q-2)}\neq0$ and $(1,q-2)$ is pure. Lemma \ref{jiben2} (i) implies $k_{1,q-2}=1$ or $2$.

\textbf{Case 2.2.1.} $k_{1,q-2}=1$.

Since $\partial(z,x_{l-1})\geq q-2$, $l\geq q-1$. Pick a vertex $x_{l-2}'\in P_{(q-1,1),(1,q-2)}(x_{l-1},w)$. The fact that $(1,q-2)$ is pure implies that $\wz{\partial}(x_{l-2}',z)=(2,q-3)$ and $l=q-1$, contrary to $q-1\notin\{2,q,p\}$.

\textbf{Case 2.2.2.} $k_{1,q-2}=2$.

Since $\partial(w,x_{l-2})<l$, one gets $h_{l-2}\neq p$ from Theorem \ref{Main1} (i). Then $h_{j}\neq p$ for $0\leq j\leq l-2$. Step \ref{step} implies $|\{i\mid h_{i}=q\}|<2$. It follows from Step \ref{(1,g-1) or (1,q-1)} that $h_{l-1}=p$ and $h_{j}=2$ for $1\leq j\leq l-2$.

Since $(q,p)\notin\mathscr{D}$ and $(1,q-2)$ is pure, by Lemma \ref{jiben4}, one has $(q-1,p)\notin\mathscr{B}\cup\mathscr{C}\cup\mathscr{D}$. In view of $w\in P_{(1,q-1),(1,q-2)}(x_{l-1},z)$, there exists a vertex $w'\in P_{(1,q-2),(1,q-1)}(x_{l-1},z)$. It follows from Theorem \ref{Main1} (iii) that $w'\neq y'$. Observe $\partial(x_{1},x)<l$. Theorem \ref{Main1} (i) and Lemma \ref{jiben2} (i) imply $\wz{\partial}(y',x_{1})\neq(1,q-2)$ and $\wz{\partial}(w',x_{1})=(1,q-2)$. Since $(1,q-2)$ is pure, by Lemma \ref{pure-g=2}, we have $\wz{\partial}(x_{l-1},x_{1})=(2,2)$ and $l=q=4$. In view of $\partial(x_{2},w')\leq2$ and Lemma \ref{(1,g-1) (1,1)=(2,g)} (i), we obtain $\partial(w',x_{2})=1$. Then $(x_{2},x_{3},w')$ is a circuit containing arcs of types $(1,1)$ and $(1,2)$, a contradiction.

\textbf{Case 2.3.} $(t,q)\in\mathscr{C}\cup\mathscr{D}$.

Observe that $(1,t-1)$ is pure and $(1,t)$ is mixed. By Lemma \ref{jiben4}, we have $k_{1,t-1}=1$ or $\Delta_{t}\simeq{\rm Cay}(\mathbb{Z}_{2t},\{1,t+1\})$. It follows from Theorem \ref{Main1} (ii) that $p_{(1,t),(1,t)}^{(1,t-1)}\neq0$. Pick a vertex $x_{l-2}'\in P_{(t,1),(1,t-1)}(x_{l-1},w)$. In view of $k_{1,t-1}=1$ or Lemma \ref{jiben3} (i), one gets $\wz{\partial}(x_{l-2}',z)=(2,t-2)$. Note that $l-1=\partial(z,x_{l-1})\geq t-1$. Hence, $l=t$. Since $(1,t-1)$ is pure, we obtain $t=p$, contrary to $(q,p)\notin\mathscr{C}\cup\mathscr{D}$.

We complete the proof of Step \ref{pre-main3 1}.

\begin{step}\label{pre-main3 2}
~{\rm Show that $(1,s-1)$ is pure if $\partial(x_{l-1},z)=2$.}\vspace{-0.3cm}
\end{step}

Suppose for the contrary that $(1,s-1)$ is mixed. Theorem \ref{Main1} (ii) implies that $p_{(1,s-1),(1,s-1)}^{(1,s-2)}\neq0$ and $(1,s-2)$ is pure. Pick vertices $w\in P_{(1,s-1),(1,1)}(x_{l-1},z)$ and $x_{l-2}'\in P_{(s-1,1),(1,s-2)}(x_{l-1},w)$.

\textbf{Case 1.} $s>3$.

By Lemma \ref{mix-(2,q-2)'} (ii), we have $\partial(z,x_{l-1})\neq s-2$. Since $(q,p)\notin\mathscr{B}$, one gets $\{s-1,s\}\neq\{q,p\}$. It follows from Lemma \ref{(1,g-1) (1,1)=(2,g)} (ii) that $\wz{\partial}(x_{l-1},z)=(2,s-1)$ and $l=s$. Since $(1,s-2)$ is pure, by Lemma \ref{pure-g=2} and Theorem \ref{Main1} (i), we obtain $\partial(x_{l-2}',z)=2$. From Lemma \ref{(1,g-1) (1,1)=(2,g)} (i), we obtain $\wz{\partial}(x_{l-2}',z)=(2,s-1)$. Since $z\in P_{(2,s-1),(s-1,2)}(x_{l-2}',x_{l-1})$, by Lemma \ref{jiben2} (iv), one has $s=2$, a contradiction.

\textbf{Case 2.} $s=3$.

Note that $l\leq4$ and $x_{l-2}'\neq z$. Lemma \ref{jiben2} (i) implies $k_{1,1}=2$. Without loss of generality, we may assume $p\neq 3$. By Lemma \ref{mix-(2,q-2)'} (ii) and Theorem \ref{Main1} (ii), we have $p=4$ or $5$.

Suppose $p=4$. Theorem \ref{Main1} (ii) implies that $(1,3)$ is pure. If $\wz{\partial}(x_{l-2}',z)=(2,2)$, by $p_{(1,3),(1,3)}^{(2,2)}\neq0$, then there exists a vertex $w'\in P_{(1,3),(1,3)}(x_{l-2}',z)$, which implies that $(x_{l-2}',w',z,w)$ is a circuit containing arcs of types $(1,1)$ and $(1,3)$, a contradiction. It follows from Lemma \ref{pure-g=2} that $\wz{\partial}(x_{l-2}',z)=(1,1)$. Since $\partial(x_{l-2}',y)\leq3$, by Theorem \ref{Main1} (i) and Lemma \ref{(1,g-1) (1,1)=(2,g)} (i), we get $\wz{\partial}(y,x_{l-2}')=(1,3)$. Then $\partial(x,x_{l-2}')=2$, contrary to $l\geq 3$.

Suppose $p=5$. By Lemma \ref{mix-(2,q-2)'} (ii) and Theorem \ref{Main1} (ii), one has $l=4$. Since $\partial(w,x)\leq3$, we obtain $\partial(y,w)=2$ from Lemma \ref{pure-g=2} and Theorem \ref{Main1} (i). In view of $\partial(w,y)\leq4$ and Lemma \ref{(1,g-1) (1,1)=(2,g)} (i), $(1,4)$ is mixed. It follows from Theorem \ref{Main1} (ii) that $(1,3)$ is pure. By Lemma \ref{pure-g=2}, we get $\wz{\partial}(x_{l-2}',z)=(2,2)$, which implies that there exists a vertex $w'\in P_{(1,3),(1,3)}(x_{l-2}',z)$. Hence, $(x_{l-2}',w',z,w)$ is a circuit containing arcs of types $(1,1)$ and $(1,3)$, a contradiction.

\begin{step}\label{pre-main3 3}
~{\rm Show that $\{h_{l-1},s\}=\{q,p\}$ if $\partial(x_{l-1},z)=2$.}\vspace{-0.3cm}
\end{step}

By Step \ref{pre-main3 2}, $(1,s-1)$ is pure. From Step \ref{pre-main3 1} and Lemma \ref{(1,g-1) (1,1)=(2,g)} (i), we get $\wz{\partial}(x_{l-1},z)=(2,s)$ and $l=s+1$, which imply $s\in\{q,p\}$. Pick a vertex $w\in P_{(1,s-1),(1,1)}(x_{l-1},z)$. Without loss of generality, we may assume $s=q$.

Suppose $h_{l-1}=2$. Observe that $\wz{\partial}(x_{l-1},z)=(2,q)=(2,l-1)$. By Lemma \ref{pure-g=2} and Theorem \ref{Main1} (i), we get $\partial(x_{l-1},y)=2$. In view of Lemma \ref{(1,g-1) (1,1)=(2,g)} (i), one has $\wz{\partial}(x_{l-1},y)=(2,q)$. It follow from Lemma \ref{jiben2} (iv) that $p=2$, a contradiction. Hence, $h_{l-1}\in\{q,p\}$.

Suppose $h_{l-1}=s$. Since $l\geq3$, one gets $w\neq x$ and $k_{1,q-1}=2$ from Lemma \ref{jiben2} (i). By Lemma \ref{jiben2} (iv) and $x_{l-1}\in P_{(q-1,1),(1,q-1)}(x,w)$, we have $\partial(x,w)=\partial(w,x)$. In view of $z\in P_{(2,l),(1,1)}(x,w)$, there exists $z'\in P_{(1,1),(2,l)}(w,x)$. Since $l\geq3$, we get $z'\neq z$ and $k_{1,1}=2$. By $\wz{\partial}(w,x)\neq(1,1)$, we obtain $p_{(1,q-1),(q-1,1)}^{(1,1)}=p_{(1,q-1),(1,1)}^{(1,q-1)}=0$. In view of Lemma \ref{pure-g=2}, Theorem \ref{Main1} (i) and Lemma \ref{(1,g-1) (1,1)=(2,g)} (i), one has $\Gamma_{1,q-1}\Gamma_{1,1}=\{\Gamma_{2,q}\}$. Since $p_{(1,q-1),(1,1)}^{(2,q)}=2$ from Lemma \ref{jiben} (i), one obtains $\wz{\partial}(x,z)=(1,1)$, a contradiction. Thus, $h_{l-1}=p$.

\begin{step}\label{(x_l-1,z) (3,l-1)1}
~{\rm For $a,b\in\{2,q,p\}$ and $a<b$, show that $p_{(1,b-1),(2,l)}^{(3,l-1)}\neq0$ if $p_{(1,a-1),(2,l)}^{(3,l-1)}=0$.}\vspace{-0.3cm}
\end{step}

Without loss of generality, we may assume $b=q$. We claim that $h_{i}=q$ for some $i\in\{0,1,\ldots,l-1\}$. Assume the contrary, namely, $h_{i}\neq q$ for each $i$. Suppose $a=2$. By Step \ref{(1,g-1) or (1,q-1)}, we may assume $h_{l-1}=2$. It follows from Steps \ref{partial neq 1}, \ref{pre-main3 1} and \ref{pre-main3 3} that $\wz{\partial}(x_{l-1},z)=(3,l-1)$, contrary to $p_{(1,1),(2,l)}^{(3,l-1)}=0$. Suppose $a=p$. By Step \ref{(1,g-1) or (1,q-1)}, we may assume
$h_{l-2}=h_{l-1}=p$. It follows from Steps \ref{partial neq 1}, \ref{pre-main3 1} and \ref{pre-main3 3} that $\partial(x_{l-1},z)=2$ and there exists a vertex $w\in P_{(1,q-1),(1,1)}(x_{l-1},z)$. Theorem \ref{Main1} (i) implies $\partial(x_{l-2},w)=2$, contrary to $\partial(w,x_{l-2})<l$. So our claim is valid.

Without loss of generality, we may assume $h_{l-1}=q$. It suffices to show that $\partial(x_{l-1},z)=3$.

Suppose $\partial(x_{l-1},z)=2$. It follows from Steps \ref{pre-main3 1}--\ref{pre-main3 3} that $(1,p-1)$ is pure and there exists a vertex $w'\in P_{(1,1),(1,p-1)}(x_{l-1},z)$. By Lemma \ref{(1,g-1) (1,1)=(2,g)} (i), we have $\wz{\partial}(x_{l-1},z)=(2,p)$ and $l=p+1$. Let $(y_{0}=z,y_{1},\ldots,y_{l-2}=w')$ be a path consisting of arcs of type $(1,p-1)$. Since $x_{l-1}\in P_{(1,1),(1,q-1)}(w',x)$, there exists $x_{l-1}'\in P_{(1,q-1),(1,1)}(w',x)$. Note that $(z=y_{0},y_{1},\ldots,y_{l-2},x_{l-1}',x)$ is a shortest path. Then $\wz{\partial}(y_{l-3},x_{l-1}')=(2,l)$. Hence, $\wz{\partial}(x,y_{1})=\wz{\partial}(x_{l-1}',z)=(3,l-1)$, contrary to $a\in\{2,p\}$. By Step \ref{partial neq 1}, we obtain $\partial(x_{l-1},z)=3$, as desired.

\vspace{3ex}

Based on the above discussion, we consider two cases, and reach a contradiction, respectively.

\textbf{Case 1.}~$p_{(2,l),(1,1)}^{(3,l-1)}\neq0$.

Pick a vertex $y_{1}\in P_{(3,l-1),(1,1)}(x,z)$. By Step \ref{(x_l-1,z) (3,l-1)1}, we may assume $p_{(1,q-1),(2,l)}^{(3,l-1)}\neq0$. Then there exist vertices $z'\in P_{(2,l),(1,q-1)}(x,y_{1})$ and $y_{1}'\in P_{(1,p-1),(1,q-1)}(x,z')$. It follows from Lemma \ref{jiben2} (i) that $k_{2,l}=2$. Observe that $x\in P_{(l,2),(2,l)}(z,z')$. By Lemma \ref{jiben2} (iv) and Lemma \ref{pure-g=2}, we get $\wz{\partial}(z',z)=(2,2)$. Lemma \ref{mix-(2,q-2)'} (ii) and Theorem \ref{Main1} (ii) imply $q=3$.

By $\wz{\partial}(z',z)=(2,2)$ and Lemma \ref{(1,g-1) (1,1)=(2,g)} (i), $(1,2)$ is mixed, which implies $p_{(1,2),(1,2)}^{(1,1)}\neq0$ from Theorem \ref{Main1} (ii). Since $\partial(y_{1}',y_{1})=2$, by Lemma \ref{jiben2} (ii), we have $|(\Gamma_{1,2})^{2}|=2$ and $k_{1,2}=2$. In view of $\wz{\partial}(z,z')=(2,2)$, $p_{(1,2),(1,2)}^{(1,1)}=1$. It follows from Lemma \ref{jiben2} (iii) that there exists a vertex $z''\in P_{(1,2),(1,2)}(y_{1}',y_{1})\setminus\{z'\}$. In view of Theorem \ref{Main1} (i), we get $z''\in\Gamma_{2,l}(x)$. Since $k_{2,l}=2$, we obtain $z''=z$, a contradiction.

\textbf{Case 2.}~$p_{(2,l),(1,1)}^{(3,l-1)}=0$.

We claim that any shortest path from $z$ to $x$ does not contain an edge. Suppose for the contrary that $h_{l-1}=2$. It follows from Steps \ref{partial neq 1}, \ref{pre-main3 1} and \ref{pre-main3 3} that $\partial(x_{l-1},z)=3$, contrary to $p_{(2,l),(1,1)}^{(3,l-1)}=0$. Thus, the claim is valid. By Step \ref{(x_l-1,z) (3,l-1)1}, we have $p_{(1,q-1),(2,l)}^{(3,l-1)}\neq0$ and $p_{(1,p-1),(2,l)}^{(3,l-1)}\neq0$. Pick a vertex $y_{l-1}\in P_{(q-1,1),(3,l-1)}(x,z)$.  It follows that there exist vertices $x'\in P_{(1,p-1),(2,l)}(y_{l-1},z)$ and $y''\in P_{(1,p-1),(1,q-1)}(x',z)$. By Lemma \ref{jiben2} (i), $k_{2,l}=2$. In view of Lemma \ref{jiben} (i), one obtains $k_{1,q-1}=2$ or $k_{1,p-1}=2$.

\textbf{Case 2.1}.~$k_{1,q-1}=2$ and $k_{1,p-1}=2$.

In view of the claim and Step \ref{(1,g-1) or (1,q-1)}, there exists a vertex $z_{1}$ such that $\wz{\partial}(z,z_{1})=(1,p-1)$ and $\partial(z_{1},y_{l-1})=l-2$. By Theorem \ref{Main1} (i), if $y'=y''$, then $x,x'\in\Gamma_{l,2}(z)$; if $y'\neq y''$, then $y',y''\in\Gamma_{l,2}(z_{1})$. In view of Lemma \ref{jiben} (i), we have $A_{1,q-1}A_{1,p-1}=2A_{2,l}$, which implies $\wz{\partial}(x',y)=(1,q-1)$. Thus, $\wz{\partial}(y_{l-1},y)=(2,l)$ and $\wz{\partial}(x,y)=(1,p-1)$, a contradiction.

\textbf{Case 2.2}.~$k_{1,q-1}=1$ or $k_{1,p-1}=1$.

Without loss of generality, we may assume $k_{1,q-1}=1$. Then $y'=y''$ and $k_{1,p-1}=2$. It follows from Theorem \ref{Main1} (i) that $\wz{\partial}(y_{l-1},y')=(2,l)$. In view of $x'\in P_{(1,p-1),(1,p-1)}(y_{l-1},y')$, one gets $p_{(1,p-1),(1,p-1)}^{(2,l)}=1$ and there exists a vertex $y_{0}\in P_{(1,p-1),(1,p-1)}(x,z)$. Since $k_{2,l}=2$, by Lemma \ref{jiben} (i) and Lemma \ref{jiben2} (ii), we have $|(\Gamma_{1,p-1})^{2}|=2$. In view of Theorem \ref{Main1} (i), we get $y_{0}\in\Gamma_{2,l}(y_{l-1})$. Then $\wz{\partial}(x',y_{0})\neq(1,p-1)$. Since $p_{(1,p-1),(1,p-1)}^{(2,l)}=1$, we obtain $y\in P_{(1,p-1),(1,p-1)}(x',z)$ and $\wz{\partial}(y_{l-1},y)\neq(2,l)$. By $\partial(y_{l-1},y)=2$ and Theorem \ref{Main1} (ii), $(1,p-1)$ is pure, which implies $l>p-2$. Then $\wz{\partial}(y_{l-1},y)=(2,p-2)$, contrary to $x\in P_{(1,q-1),(1,q-1)}(y_{l-1},y)$.

Thus, we finish the proof of Proposition \ref{jiegou} for the case B.

\section{Subdigraphs}

In this section, we focus on the existence of some special subdigraphs of commutative quasi-thin weakly distance-regular digraphs.

Let $F$ be a nonempty subset of $R$ and $x\in V\Gamma$. Set $F(x):=\{y\in V\Gamma\mid(x,y)\in \cup_{f\in F}f\}$, and $F_{q_{1},q_{2},\ldots,q_{l}}(x)$ is a collection of vertices $y$ satisfying each arc in one of paths from $x$ to $y$ is of type $(1,q_{1}-1),(1,q_{2}-1),\ldots,(1,q_{l-1}-1)$ or $(1,q_{l}-1)$. If $\Gamma_{\wz{i}^{*}}\Gamma_{\wz{j}}\subseteq F$ for any $\Gamma_{\wz{i}},\Gamma_{\wz{j}}\in F$, we say that $F$ is {\em closed}. Let $\langle F\rangle$ be the minimum closed subset containing $F$. We write $\langle\Gamma_{1,q-1}\rangle$ instead of $\langle\{\Gamma_{1,q-1}\}\rangle$.

\begin{prop}\label{C_q,h pre1}
If $C_{q,h}$ exists, then $\Delta_{q,h}\simeq{\rm Cay}(\mathbb{Z}_{q}\times\mathbb{Z}_{4},\{(1,0),(0,1),(1,2)\})$ for $q>2$ and $q\neq h$.
\end{prop}
\textbf{Proof.}~For fixed $x\in V\Gamma$, by Lemma \ref{jiben4}, there exists an isomorphism $\tau$ from ${\rm Cay}(\mathbb{Z}_{2q},\{1,q+1\})$ to $\Delta_{q}(x)$. Write $\tau(a)=(a,0)$ for each $a\in\mathbb{Z}_{2q}$. Suppose that there exists a vertex $(s,0)\in\Gamma_{1,h-1}(0,0)$. From Lemma \ref{jiben3} (i), we have $s=q$. Since $(1,0)\in P_{(1,q-1),(q-1,1)}((0,0),(q,0))$ , by Lemma \ref{jiben2} (iv), we get $h=2$, contrary to $h\geq3$. Hence, $\Gamma_{1,h-1}\notin\langle\Gamma_{1,q-1}\rangle$. In view of Lemma \ref{jiben4}, one obtains $k_{1,h-1}=1$. Since $C_{q,h}$ exists, $V\Delta_{q,h}(x)$ has a partition $F_{q}(x)\dot{\cup}F_{q}(x')$. It follows that $\sigma:F_{q}(x)\rightarrow F_{q}(x')$, $y\mapsto y'$ is an isomorphism from $\Delta_{q}(x)$ to $\Delta_{q}(x')$, where $y'\in \Gamma_{1,h-1}(y)$. Write $\sigma(a,0)=(a,1)$ for each $a$. Since $C_{q,h}$ exists again, $((a,1),(a+q,0))\in\Gamma_{1,h-1}$. The desired result holds.$\qed$

\begin{prop}\label{pure}
Let $q\geq3$. If $k_{1,q-1}=2$ and $(1,q-1)$ is pure, then $\Delta_{q}\simeq{\rm Cay}(\mathbb{Z}_{2q},\{1,q+1\})$.
\end{prop}
\textbf{Proof.} Suppose not. By Lemma \ref{jiben3}, there exists an isomorphism $\tau$ from ${\rm Cay}(\mathbb{Z}_{q}\times\mathbb{Z}_{q},\{(1,0),(0,1)\})$ to $\Delta_{q}$.

By Lemma \ref{jiben4} and Proposition \ref{jiegou}, C4, C5 or C6 holds, which implies that $K\subseteq\{(1,1),(1,q-1),(1,q)\}$. If $(1,q)\in\wz{\partial}(\Gamma)$, then $(1,q)$ is mixed, which implies $p_{(1,q),(1,q)}^{(1,q-1)}\neq0$ and $k_{1,q}=2$ by Lemma \ref{jiben} (i), Lemma \ref{jiben2} (i) and Theorem \ref{Main1} (ii).

\begin{stepp}\label{pure-step 1}
~{\rm Show that $\Gamma_{1,q-1}\Gamma_{1,1}=\{\Gamma_{2,q}\}$ if $(1,1)\in\wz{\partial}(\Gamma)$.}\vspace{-0.3cm}
\end{stepp}

Suppose $p_{(1,q-1),(q-1,1)}^{(1,1)}\neq0$. Note that $\wz{\partial}_{\Gamma}(\tau(a,b),\tau(a+1,b-1))=(1,1)$. By Lemma \ref{jiben2} (i), $k_{1,1}=2$. Observe that $\tau(1,0)\in P_{(1,1),(1,q-1)}(\tau(0,1),\tau(2,0))$ and $(\tau(0,1),\tau(2,0))\notin\Gamma_{1,q-1}\cup\Gamma_{1,1}$. In view of Theorem \ref{Main1} (i) and Lemma \ref{(1,g-1) (1,1)=(2,g)} (i), we get $(\tau(0,1),\tau(2,0))\in\Gamma_{2,q}$, contrary to the fact that $(\tau(2,0),\tau(3,0),\ldots,\tau(0,0),\tau(0,1))$ is a path of length $q-1$. Thus, $p_{(1,q-1),(q-1,1)}^{(1,1)}=0$. It follows that $\Gamma_{1,q-1}\Gamma_{1,1}=\{\Gamma_{2,q}\}$.

\begin{stepp}\label{pure-step 2}
~{\rm Show that $\Gamma_{1,q}\Gamma_{1,1}=\{\Gamma_{2,q+1}\}$ if $(1,1),(1,q)\in\wz{\partial}(\Gamma)$.}\vspace{-0.3cm}
\end{stepp}

Let $x,y,z,w$ be vertices such that $\wz{\partial}(x,y)=(1,q)$, $\wz{\partial}(y,z)=(1,1)$ and $w\in P_{(q,1),(1,q-1)}(x,y)$. By Step \ref{pure-step 1}, $\wz{\partial}(w,z)=(2,q)$. Since $k_{1,q-1}=2$, from Lemma \ref{jiben} (i) and Lemma \ref{jiben2} (i), we obtain $k_{2,q}=2$. Suppose $\partial(x,z)=1$. In view of Lemma \ref{pure-g=2} and Theorem \ref{Main1} (i), one has $\wz{\partial}(x,z)=(1,q)$. Note that $x\in P_{(1,q),(1,q)}(w,z)$ and $p_{(1,q),(1,q)}^{(1,q-1)}\neq0$. By Lemma \ref{jiben2} (ii), we get $(\Gamma_{1,q})^{2}=\{\Gamma_{1,q-1},\Gamma_{2,q}\}$. Since $k_{1,q}=2$, from Lemma \ref{jiben} (i) and (v), we obtain $p_{(1,q),(1,q)}^{(1,q-1)}=1$. In view of Lemma \ref{jiben2} (iii), we have $p_{(1,q),(1,q)}^{(2,q)}=2$, which implies $k_{2,q}=1$, a contradiction. Then $\partial(x,z)=2$. Since $\wz{\partial}(w,z)=(2,q)$, by Lemma \ref{jiben2} (iv), we have $\partial(z,x)\neq q$. In view of Lemma \ref{mix-(2,q-2)'} (ii), $\partial(z,x)=q+1$. Thus, $\Gamma_{1,q}\Gamma_{1,1}=\{\Gamma_{2,q+1}\}$.

\begin{stepp}\label{pure-step 3}
~{\rm Show that $(A_{1,q-1})^{2}=A_{2,q-2}+2A_{2,2q-2}$.}\vspace{-0.3cm}
\end{stepp}

In view of Lemma \ref{jiben3} (ii) and Theorem \ref{Main1} (ii), we have $(A_{1,q-1})^{2}=A_{2,q-2}+p_{(1,q-1),(1,q-1)}^{(2,t)}A_{2,t}$ with $t\neq q-2$. By Lemma \ref{jiben2} (iii), one gets $p_{(1,q-1),(1,q-1)}^{(2,t)}=2$, which implies $k_{2,t}=1$ from Lemma \ref{jiben} (i) and (v). Let $x,y,y',z$ be vertices such that $\wz{\partial}(x,z)=(2,t)$ and $P_{(1,q-1),(1,q-1)}(x,z)=\{y,y'\}$.

We claim that $\partial(x,x_{1})=3$ for any path $(z=x_{0},x_{1},\ldots,x_{t}=x)$. Assume the contrary, namely, $\partial(x,x_{1})=1$ or $2$.

\textbf{Case 1.} $\partial(x,x_{1})=1$.

Since $x_{1}\notin\{y,y'\}$, we have $\wz{\partial}(x,x_{1})=(1,1)$ or $(1,q)$. If $\wz{\partial}(x,x_{1})=(1,1)$, by Step \ref{pure-step 1}, then $\wz{\partial}(x_{1},y)=(2,q)$, contrary to $q>2$; if $\wz{\partial}(x,x_{1})=(1,q)$, by $p_{(1,q),(1,q)}^{(1,q-1)}\neq0$, then $y$ or $y'\in\Gamma_{1,q}(x_{1})$, which implies that $(y,z,x_{1})$ or $(y',z,x_{1})$ is a circuit, contrary to $q>2$.

\textbf{Case 2.} $\partial(x,x_{1})=2$.

Pick a vertex $w\in P_{(1,h-1),(1,l-1)}(x,x_{1})$. Suppose $h=q$. Then $w\in\{y,y'\}$. Since $(1,q-1)$ is pure, $\wz{\partial}(w,x_{1})\neq(1,1)$. In view of Theorem \ref{Main1} (i) and (ii), we have $\wz{\partial}(z,x_{1})=(1,1)$, and $y$ or $y'\in\Gamma_{q-1,1}(x_{1})$, which imply $p_{(1,q-1),(1,1)}^{(1,q-1)}\neq0$, contrary to Step \ref{pure-step 1}. Thus, $h\neq q$ and $l\neq q$.

Suppose $h=l=2$. Lemma \ref{jiben2} (i) implies $k_{1,1}=2$. By Step \ref{pure-step 1}, $y,y'\in\Gamma_{2,q}(w)$. It follows from Lemma \ref{jiben} (i) that $p_{(1,1),(1,q-1)}^{(2,q)}=2$ and $y,y'\in\Gamma_{1,q-1}(x_{1})$. Since $(1,q-1)$ is pure, we get $q=3$ and $\wz{\partial}(z,x_{1})=(1,2)$. Observe that $y,y'\in P_{(1,2),(1,2)}(x_{1},z)$, contrary to $p_{(1,2),(1,2)}^{(2,1)}=1$.

Suppose $h=q+1$ or $l=q+1$. By $p_{(1,q),(1,q)}^{(1,q-1)}\neq0$, we may assume that $h=q+1$ and $\wz{\partial}(w,y)=(1,q)$. Since $\partial(y,x_{1})\leq2$, one gets $l=q+1$ from Step \ref{pure-step 2}. In view of $\wz{\partial}(x,x_{1})=(2,t-1)$ and Lemma \ref{jiben2} (ii), one has $(\Gamma_{1,q})^{2}=\{\Gamma_{1,q-1},\Gamma_{2,t-1}\}$. Since $k_{1,q-1}=2$, by Lemma \ref{jiben} (i) and (v), we obtain $p_{(1,q),(1,q)}^{(1,q-1)}=1$. By Lemma \ref{jiben2} (iii), we get $p_{(1,q),(1,q)}^{(2,t-1)}=2$, which implies $k_{2,t-1}=1$. Since $k_{1,q-1}=k_{1,q}=2$ and $k_{2,t}=1$, from Lemma \ref{jiben} (i), one has $\wz{\partial}(z,x_{1})=(1,1)$. In view of Step \ref{pure-step 1}, $\wz{\partial}(y,x_{1})=(2,q)$. Since $w\in P_{(q,1),(1,q)}(y,x_{1})$, from Lemma \ref{jiben2} (iv), we get $q=2$, a contradiction.

Thus, our claim is valid.

Suppose that the path $(x_{0},x_{1},\ldots,x_{t})$ contains arcs of different types. Without loss of generality, we may assume $\wz{\partial}(z,x_{1})=(1,u-1)$ and $\wz{\partial}(x_{1},x_{2})=(1,v-1)$ with $u\neq v$. Pick a vertex $x_{1}'\in P_{(1,v-1),(1,u-1)}(z,x_{2})$. By the claim, we get $\wz{\partial}(x,x_{1})=\wz{\partial}(x,x_{1}')=(3,t-1)$. It follows from Lemma \ref{jiben} (iv) that $k_{2,t}\geq2$, a contradiction. Then the path $(x_{0},x_{1},\ldots,x_{t})$ consists of arcs of the same type.

Suppose $\wz{\partial}(x_{i},x_{i+1})=(1,1)$ for $0\leq i\leq t-1$. By Lemma \ref{pure-g=2}, $t=2$. In view of Step \ref{pure-step 1}, we get $\wz{\partial}(y,x_{1})=(2,q)$. Since $(x_{1},x_{2}=x,y)$ is a path, one has $q\leq2$, a contradiction.

Suppose $\wz{\partial}(x_{i},x_{i+1})=(1,q)$ for $0\leq i\leq t-1$. Then $\partial(z,x_{2})=2$. In view of $p_{(1,q),(1,q)}^{(1,q-1)}\neq0$ and Lemma \ref{mix-(2,q-2)'} (i), we have $\partial(x_{2},z)>q-1$, which implies $t\geq3$. Since $k_{1,q-1}=2$ and $|(\Gamma_{1,q})^{2}|=2$ from Lemma \ref{jiben2} (ii), one gets $p_{(1,q),(1,q)}^{(1,q-1)}=1$ by Lemma \ref{jiben} (i) and (v). In view of Lemma \ref{jiben2} (iii), there exists a vertex $x_{1}''\in P_{(1,q),(1,q)}(z,x_{2})$ such that $\wz{\partial}(x_{1}'',x_{3})=(1,q-1)$, a contradiction.

Hence, $\wz{\partial}(x_{i},x_{i+1})=(1,q-1)$ for $0\leq i\leq t-1$. Since $\Delta_{q}\simeq{\rm Cay}(\mathbb{Z}_{q}\times\mathbb{Z}_{q},\{(1,0),(0,1)\})$, we have $t=2q-2$.

\vspace{3ex}

In the following, we reach a contradiction based on the above discussion.

Suppose $q>3$. Note that $\wz{\partial}_{\Gamma}(\tau(a,b),\tau(a+1,b+1))=(2,2q-2)$. Since
\begin{eqnarray}
&&(\tau(1,1),\tau(2,1),\tau(3,1),\ldots,\tau(-1,1),\tau(0,1),\tau(0,2),\ldots,\tau(0,0)),\nonumber\\
&&(\tau(1,1),\tau(2,1),\tau(2,2),\ldots,\tau(2,-1),\tau(2,0),\tau(3,0),\ldots,\tau(0,0))\nonumber
\end{eqnarray}
are two shortest paths, we get $\tau(3,1),\tau(2,2)\in\Gamma_{4,2q-4}(\tau(0,0))$. But $\tau(1,1)\in P_{(2,2q-2),(2,q-2)}(\tau(0,0),\tau(3,1))$ and $P_{(2,2q-2),(2,q-2)}(\tau(0,0),\tau(2,2))=\emptyset$, a contradiction. In the following, we consider $q=3$.

\textbf{Case 1}.~$(1,1)\in\wz{\partial}(\Gamma)$.

By Step \ref{pure-step 3} and Lemma \ref{jiben} (i), we have $k_{2,4}=1$. From Step \ref{pure-step 2}, one gets $(1,3)\notin\wz{\partial}(\Gamma)$. Since the valency of $\Gamma$ is more than $3$, by Lemma \ref{jiben2} (i), one has $k_{1,1}=2$. Let $x,y,z,z'$ be distinct vertices such that $\wz{\partial}(x,y)=(1,2)$ and $\wz{\partial}(y,z)=\wz{\partial}(y,z')=(1,1)$. By Step \ref{pure-step 1}, we obtain $z,z'\in\Gamma_{2,3}(x)$. In view of Lemma \ref{jiben} (i), one has $p_{(1,2),(1,1)}^{(2,3)}=2$, which implies that there exists a vertex $y'$ such that $\wz{\partial}(x,y')=(1,2)$ and $\wz{\partial}(y',z)=\wz{\partial}(y',z')=(1,1)$ with $y'\neq y$. Hence, $(y,z,y',z')$ is an undirected circuit of length $4$. By Lemma \ref{pure-g=2}, we get $\wz{\partial}(y,y')=(2,2)$ and $p_{(1,1),(1,1)}^{(2,2)}=2$. From Lemma \ref{jiben} (i) and (v), $k_{2,2}=1$. Since $x\in P_{(2,1),(1,2)}(y,y')$, we have $p_{(2,1),(1,2)}^{(2,2)}=2$, contrary to $\Delta_{3}\simeq{\rm Cay}(\mathbb{Z}_{3}\times\mathbb{Z}_{3},\{(1,0),(0,1)\})$.

\textbf{Case 2}.~$(1,1)\notin\wz{\partial}(\Gamma)$.

Note that $(1,3)\in\wz{\partial}(\Gamma)$. Pick a vertex $w\in P_{(1,3),(1,3)}(\tau(0,0),\tau(0,1))$. By Lemma \ref{jiben2} (ii), we have $|(\Gamma_{1,3})^{2}|=1$ or $2$.

\textbf{Case 2.1}.~$|(\Gamma_{1,3})^{2}|=1$.

Since $k_{1,2}=2$, by Lemma \ref{jiben} (i), we have $p_{(1,3),(1,3)}^{(1,2)}=2$ and $\wz{\partial}_{\Gamma}(w,\tau(1,0))=(1,3)$. Pick a vertex $x'\in P_{(1,3),(1,3)}(\tau(0,0),\tau(1,0))$ with $x'\neq w$. Observe $x'\in P_{(1,3),(1,3)}(\tau(0,0),\tau(0,1))$. Since $w,x'\in P_{(3,1),(1,3)}(\tau(0,1),\tau(1,0))$, from Lemma \ref{jiben} (i) and (v), we obtain $k_{\wz{\partial}_{\Gamma}(\tau(0,1),\tau(1,0))}=1$ and $|P_{(1,2),(2,1)}(\tau(0,1),\tau(1,0))|=2$, contrary to $\Delta_{3}\simeq{\rm Cay}(\mathbb{Z}_{3}\times\mathbb{Z}_{3},\{(1,0),(0,1)\})$.

\textbf{Case 2.2}.~$|(\Gamma_{1,3})^{2}|=2$.

Since $|(\Gamma_{1,3})^{2}|=2$, $(w,\tau(1,0))\notin\Gamma_{1,3}$. It follows that $P_{(1,2),(1,3)}(w,\tau(1,1))=\{\tau(0,1)\}$. By Lemma \ref{jiben} (i) and Lemma \ref{jiben2} (i),(ii), we have $|\Gamma_{1,2}\Gamma_{1,3}|=2$. In view of Theorem \ref{Main1} (i), one obtains $\wz{\partial}_{\Gamma}(w,\tau(0,2))=(2,2)$ and $\partial_{\Gamma}(w,\tau(1,1))=2$. By Step \ref{pure-step 3}, we get $p_{(1,2),(1,2)}^{(2,4)}=2$. Hence, $\partial_{\Gamma}(\tau(1,1),w)=3$ or $5$.

\textbf{Case 2.2.1}.~$\partial_{\Gamma}(\tau(1,1),w)=3$.

Pick a path $(\tau(1,1),z_{1},z_{2},w)$. Suppose that $(z_{2},w)\in\Gamma_{1,3}$. The fact that $\partial_{\Gamma}(\tau(1,1),\tau(0,0))=4$ implies $z_{2}\neq \tau(0,0)$. Since $|(\Gamma_{1,3})^{2}|=2$, from Theorem \ref{Main1} (iii) and Lemma \ref{mix-(2,q-2)'} (i), we get $(z_{2},\tau(0,1))\in\Gamma_{2,3}$, which implies $(\Gamma_{1,3})^{2}=\{\Gamma_{1,2},\Gamma_{2,3}\}$. Since $k_{1,2}=2$, by Lemma \ref{jiben} (i) and (v), we obtain $p_{(1,3),(1,3)}^{(1,2)}=1$. In view of Lemma \ref{jiben2} (iii), one has $p_{(1,3),(1,3)}^{(2,3)}=2$ and $\wz{\partial}_{\Gamma}(\tau(0,1),\tau(1,1))=(1,3)$, a contradiction.

Observe that the path $(\tau(1,1),z_{1},z_{2},w)$ consists of arcs of type $(1,2)$. Since $(\tau(0,1),\tau(1,2)),(\tau(1,1),\tau(2,2))\in\Gamma_{2,4}$, we have $z_{1}=\tau(2,1)$ and $z_{2}=\tau(0,1)$, a contradiction.

\textbf{Case 2.2.2}.~$\partial_{\Gamma}(\tau(1,1),w)=5$.

By $\wz{\partial}_{\Gamma}(w,\tau(0,2))=(2,2)$ and Lemma \ref{jiben2} (ii), $\Gamma_{1,2}\Gamma_{1,3}=\{\Gamma_{2,2},\Gamma_{2,5}\}$. Then $\tau(2,0)\in P_{(2,1),(2,5)}(\tau(0,0),w)$. Since $(\tau(1,1),\tau(1,2),\tau(2,2),\tau(2,0),\tau(0,0),w)$ and $(\tau(1,1),\tau(2,1),\tau(2,2),\tau(2,0),\tau(0,0),w)$ are two shortest paths, $\wz{\partial}_{\Gamma}(w,\tau(1,2))=\wz{\partial}_{\Gamma}(w,\tau(2,1))=(3,4)$. It follows from Step \ref{pure-step 3} and Lemma \ref{jiben} (i) that $k_{2,4}=1$. Since $\tau(0,1)\in P_{(1,3),(2,4)}(w,\tau(1,2))$, we obtain $(w,\tau(1,0))\in\Gamma_{1,3}$, contrary to $|(\Gamma_{1,3})^{2}|=2$.

This completes the proof of the proposition.$\qed$

\begin{prop}\label{Delta(q,q+1)}
Let $q>2$, $k_{1,q-1}=2$ and $(1,q-1)$ be pure. The following hold: \vspace{-0.3cm}
\begin{itemize}
\item [{\rm(i)}] If $(1,q)$ is mixed, then $\Delta_{q,q+1}\simeq{\rm Cay}(\mathbb{Z}_{4q},\{1,2,2q+1,2q+2\})$.\vspace{-0.3cm}

\item [{\rm(ii)}] If $k_{1,1}=2$, then $\Delta_{2,q}\simeq{\rm Cay}(\mathbb{Z}_{q}\times\mathbb{Z}_{4},\{(1,0),(1,2),(0,1),(0,3)\})$ for $q\neq4$.
\end{itemize}
\end{prop}
\textbf{Proof.}~Assume that $l=q+1$ and $(1,q)$ is mixed, or $l=2$ and $k_{1,1}=2$. In view of Theorem \ref{Main1} (ii), Lemma \ref{jiben} (i) and Lemma \ref{jiben2} (i), we have $k_{1,l-1}=2$. By Proposition \ref{pure}, there exists an isomorphism $\tau$ from ${\rm Cay}(\mathbb{Z}_{2q},\{1,q+1\})$ to $\Delta_{q}(x)$ for fixed $x\in V\Gamma$. Write $\tau(a)=(a,0)$ for any $a$. Suppose that there exists a vertex $(s,0)\in\Gamma_{1,l-1}(0,0)$. By Lemma \ref{jiben3} (i), we have $s=q$. Since $(1,0),(q+1,0)\in P_{(1,q-1),(q-1,1)}((0,0),(q,0))$, from Lemma \ref{jiben2} (iv), one gets $l=2$. In view of Lemma \ref{jiben} (i) and (v), we obtain  $k_{1,1}=1$, a contradiction. Hence, $\Gamma_{1,l-1}\notin\langle\Gamma_{1,q-1}\rangle$.

If $l=q+1$, by Lemma \ref{Gamma 2,q} (ii), then $(A_{1,q})^{2}=2A_{1,q-1}$; if $l=2$, by Lemma \ref{Gamma 2,q} (i) and Lemma \ref{jiben} (i), then $A_{1,q-1}A_{1,1}=2A_{2,q}$. Then $V\Delta_{l,q}(x)$ has a partition $F_{q}(x)\dot{\cup}F_{q}(x')$. Let $\sigma$ be an isomorphism from $\Delta_{q}(x)$ to $\Delta_{q}(x')$ such that $\sigma(0,0)\in\Gamma_{1,l-1}(0,0)$. Write $\sigma(a,0)=(a,1)$ for each $a$. Suppose $l=q+1$. Since $(A_{1,q})^{2}=2A_{1,q-1}$, we have  $(a,1),(a+q,1)\in\Gamma_{1,q}(a,0)$ and $(a+1,0),(a+q+1,0)\in\Gamma_{1,q}(a,1)$, which imply that (i) holds. Suppose $l=2$. Since $A_{1,q-1}A_{1,1}=2A_{2,q}$, one gets $(a,1),(a+q,1)\in\Gamma_{1,1}(a,0)$. If $q=4$, by Lemma \ref{pure-g=2}, then $(4,0),(2,0),(6,0)\in\Gamma_{2,2}(0,0)$ since $(1,3)$ is pure, contrary to Lemma \ref{jiben2} (i). Thus, (ii) holds.$\qed$

\begin{prop}\label{Delta(2,q) k 1,q-1=1}
Suppose that {\rm C6} holds. If $k_{1,1}=2$ and $k_{1,q-1}=1$, then $\Gamma_{1,q}\notin\langle\{\Gamma_{1,1},\Gamma_{1,q-1}\}\rangle$ and $\Delta_{2,q}\simeq{\rm Cay}(\mathbb{Z}_{q}\times\mathbb{Z}_{n},\{(1,0),(0,1),(0,-1)\}$ with $n\leq q-(1+(-1)^{q})/2$.
\end{prop}
\textbf{Proof.}~Since $(1,q)$ is mixed, from Theorem \ref{Main1} (ii), $(1,q-1)$ is pure. By Lemma \ref{pure-g=2}, we get $\Gamma_{1,q-1}\notin\langle\Gamma_{1,1}\rangle$. For fixed $x_{0}\in V\Gamma$, $V\Delta_{2,q}(x_{0})$ has a partition $\dot{\cup}_{i=0}^{m-1}F_{2}(x_{i})$ with $m>1$. Let $\tau$ be an isomorphism from $\textrm{Cay}(\mathbb{Z}_{n},\{1,n-1\})$ to $\Delta_{2}(x_{0})$.  Write $\tau(a)=(0,a)$ for each $a$. Since $k_{1,q-1}=1$,  $\sigma_{j}:F_{2}(x_{j})\rightarrow F_{2}(x_{j+1})$, $y_{j}\mapsto y_{j+1}$ is an isomorphism from $\Delta_{2}(x_{j})$ to $\Delta_{2}(x_{j+1})$, where $y_{j+1}\in\Gamma_{1,q-1}(y_{j})$ for $0\leq j\leq m-2$. Write $\sigma_{j}(j,a)=(j+1,a)$.

Assume that $(s,t)\in\Gamma_{1,q-1}(m-1,0)$. Since $k_{1,q-1}=1$, we have $s=0$. It follows from Lemma \ref{jiben} (i) and Lemma \ref{pure-g=2} that $t=0$, or $2\mid n$ and $t=n/2$.

Suppose $2\mid n$ and $t=n/2$. Since $(1,q-1)$ is pure and $k_{1,q-1}=1$, from Lemma \ref{jiben} (i), we get $\wz{\partial}_{\Gamma}((0,0),(0,n/2))=(m,q-m)$, which implies $q=2m$ from Lemma \ref{pure-g=2}. Hence, $q\leq n$ and $(0,0)\in\Gamma_{1,q-1}(m-1,n/2)$. Since $\{i\mid(1,i-1)\in\wz{\partial}(\Gamma)\}=\{2,q,q+1\}$, one has $(0,m),(0,-m)\in\Gamma_{m,m}(0,0)$. Since $k_{m,m}\leq2$ by Lemma \ref{jiben2} (i), we obtain $m=n/2$ and $n=q$. Hence, $((0,0),(1,0),\ldots,(m-1,0),(0,n/2),(0,n/2-1),\ldots,(0,1))$ is a circuit of length $q$ containing arcs of types $(1,1)$ and $(1,q-1)$, contrary to the fact that $(1,q-1)$ is pure. Then $t=0$ and $m=q$. Since $(1,q-1)$ is pure and $k_{1,q-1}=1$, one has $((m-1,a),(0,a))\in\Gamma_{1,q-1}$ for each $a$. Thus, $\Delta_{2,q}\simeq{\rm Cay}(\mathbb{Z}_{q}\times\mathbb{Z}_{n},\{(1,0),(0,1),(0,-1)\}$.

Since $(1,q)$ is mixed, we have $p_{(1,q),(1,q)}^{(1,q-1)}=k_{1,q}$ from Theorem \ref{Main1} (ii) and Lemma \ref{jiben} (ii). We prove $n\leq q-(1+(-1)^{q})/2$ by way of contradiction. Assume that $n>q-(1+(-1)^{q})/2$. Suppose that $q$ is even. Since $(1,q-1)$ is pure and $k_{1,q-1}=1$, by Lemma \ref{jiben} (i), we get $\wz{\partial}_{\Gamma}((0,0),(q/2,0))=(q/2,q/2)$ and $k_{q/2,q/2}=1$. Observe $\wz{\partial}_{\Gamma}((0,0),(0,q/2))=(q/2,q/2)$, a contradiction. Suppose that $q$ is odd. Pick a vertex $x\in P_{(1,q),(1,q)}(((q-1)/2,0),((q+1)/2,0))$. Note that  $x,(0,(q+1)/2)\in\Gamma_{(q+1)/2,(q+1)/2}(0,0)$. Since $(x,((q+1)/2,0),((q+3)/2,0),\ldots,(0,0))$ is a path containing arcs of types $(1,q-1)$ and $(1,q)$, there exists a path $((0,(q+1)/2)=x_{0},x_{1},\ldots,x_{(q+1)/2}=(0,0))$ containing arcs of types $(1,q-1)$ and $(1,q)$. Then $((0,0),(0,1),\ldots,(0,(q+1)/2)=x_{0},x_{1},\ldots,x_{(q-1)/2})$ is a circuit of length $q+1$ containing arcs of types $(1,1),(1,q-1)$ and $(1,q)$, contrary to Lemma \ref{mix-(2,q-2)'} (ii).

Suppose that $(h,l)\in\Gamma_{1,q}(0,0)$ for some $h\in\{0,1,\ldots,q-1\}$ and $l\in\mathbb{Z}_{n}$. By Lemma \ref{pure-g=2}, $h\neq0$. Without loss of generality, we may assume $2\hat{l}\leq n$. The fact that $p_{(1,q),(1,q)}^{(1,q-1)}=k_{1,q}$ implies $\wz{\partial}_{\Gamma}((h,l),(1,0))=(1,q)$. Since $((0,0),(h,l),(h+1,l),\ldots,(0,l),(0,l-1),\ldots,(0,1))$ and $((1,0),(2,0),\ldots,(h,0),(h,1),\ldots,(h,l))$ are two circuits, one has $q-h+\hat{l}+1\geq q+1$ and $h+\hat{l}\geq q+1$. Hence, $q+1\leq2\hat{l}\leq n$, contrary to $n\leq q-(1+(-1)^{q})/2$. Thus, $\Gamma_{1,q}\notin\langle\{\Gamma_{1,1},\Gamma_{1,q-1}\}\rangle$.$\qed$

\section{Proof of Theorem \ref{Main}}

For any nonempty subset $F$ of $R$ with $F=\langle F\rangle$, let
\begin{eqnarray}
V\Gamma/F:=\{F(x)\mid x\in V\Gamma\}\quad\textrm{and}\quad \Gamma_{\wz{i}}^{F}:=\{(F(x),F(y))\mid y\in F\Gamma_{\wz{i}}F(x)\}.\nonumber
\end{eqnarray}
The digraph $(V\Gamma/F,\ \cup_{(1,s)\in\wz{\partial}(\Gamma)}\Gamma_{1,s}^{F})$ is said to be
the \emph{quotient digraph} of $\Gamma$ over $F$, denoted by $\Gamma/F$.

In the following, we divide the proof of Theorem \ref{Main} into four subsections according to separate assumptions based on Proposition \ref{jiegou}.

\subsection{The cases C1, C2 and C3}

By Lemma \ref{jiben4}, $k_{1,q-1}=2$. If $(1,1)\in\wz{\partial}(\Gamma)$, by Lemma \ref{C_q,h holds k_1,1} (i), then $k_{1,1}=1$; if $(1,q)\in\wz{\partial}(\Gamma)$, then $(1,q)$ is mixed, which imply $p_{(1,q),(1,q)}^{(1,q-1)}\neq0$ and $k_{1,q}=2$ from Theorem \ref{Main1} (ii) and Lemma \ref{jiben4}.

\textbf{Case 1.} $(1,q)\notin\wz{\partial}(\Gamma)$.

Note that C1 holds. Since $C_{q,3}$ exists, from Lemma \ref{h1}, $(1,2)$ is mixed. By Lemma \ref{jiben4}, we have $k_{1,2}=1$, which implies $p_{(1,2),(1,2)}^{(1,1)}=1$ from Theorem \ref{Main1} (ii). In view of Proposition \ref{C_q,h pre1}, $\Gamma$ is isomorphic to one of the digraphs in Theorem \ref{Main} {\rm(iv)} for $i=0$.

\textbf{Case 2.} $(1,q)\in\wz{\partial}(\Gamma)$.

Note that C2 or C3 holds. Assume that $h=4$ or $3$. Since $C_{q,h}$ exists, by Proposition \ref{C_q,h pre1}, there exists an isomorphism $\tau$ from $\textrm{Cay}(\mathbb{Z}_{q}\times\mathbb{Z}_{4},\{(1,0),(0,1),(1,2)\})$ to $\Delta_{q,h}(x)$ for fixed $x\in V\Gamma$. Write $\tau(a,b)=(a,b,0)$ for each $(a,b)$. Suppose that there exists $(c,d,0)$ such that $\wz{\partial}_{\Gamma}((0,0,0),(c,d,0))=(1,q)$. Since $p_{(1,q),(1,q)}^{(1,q-1)}\neq0$ and $(1,q-1)$ is pure from Lemma \ref{jiben4}, we get $d\in\{1,3\}$ and $c\neq0$. Observe that $((0,0,0),(c,d,0),(c+1,3,0),(c+2,3,0),\ldots,(0,3,0))$ is a circuit of length $q-\hat{c}+2$ containing arcs of types $(1,q)$ and $(1,h-1)$, contrary to Lemma \ref{mix-(2,q-2)'} (ii). Hence, $\Gamma_{1,q}\notin\langle\{\Gamma_{1,q-1},\Gamma_{1,h-1}\}\rangle$.

By Lemma \ref{jiben4} and Lemma \ref{Gamma 2,q} (ii), we have $(A_{1,q})^{2}=2A_{1,q-1}$, which implies that $V\Delta_{q,q+1,h}$ has a partition $F_{q,h}(x)\dot{\cup}F_{q,h}(x')$. Let $\sigma$ be an isomorphism from $\Delta_{q,h}(x)$ to $\Delta_{q,h}(x')$ such that $\sigma(0,0,0)\in\Gamma_{1,q}(0,0,0)$. Write $\sigma(a,b,0)=(a,b,1)$ for each $(a,b)$. Since $(0,0,1)\in P_{(1,q),(1,h-1)}((0,0,0),(0,1,1))$ and $k_{1,h-1}=1$, we get $(0,1,1)\in\Gamma_{1,q}(0,1,0)$. Similarly, $(0,b,1)\in\Gamma_{1,q}(0,b,0)$ for each $b$. The fact that $(A_{1,q})^{2}=2A_{1,q-1}$ implies that $(a,b,1),(a,b+2,1)\in\Gamma_{1,q}(a,b,0)$ and $(a+1,b,0),(a+1,b+2,0)\in\Gamma_{1,q}(a,b,1)$ for each $(a,b)$. Thus, $\Delta_{q,q+1,h}\simeq{\rm Cay}(\mathbb{Z}_{2q}\times\mathbb{Z}_{4},\{(2,0),(2,2),(1,0),(1,2),(0,1)\})$.

If C2 holds, then $\Gamma$ is isomorphic to one of the digraphs in Theorem \ref{Main} (vii) for $i=1$. Suppose that C3 holds. Since $C_{q,3}$ exists, from Lemma \ref{h1}, $(1,2)$ is mixed. By Lemma \ref{jiben4}, we have $k_{1,2}=1$, which implies $p_{(1,2),(1,2)}^{(1,1)}=1$ from Theorem \ref{Main1} (ii). Hence, $\Gamma$ is isomorphic to one of the digraphs in Theorem \ref{Main} (vii) for $i=0$.

We complete the proof of the main theorem for the cases C1, C2 and C3.

\subsection{The case C4}

Since the valency of $\Gamma$ is more than $3$, from Lemma \ref{jiben2} (i), we have $k_{1,1}=k_{1,q-1}=2$. By Proposition \ref{Delta(q,q+1)} (ii), $\Gamma$ is isomorphic to one of the digraphs in Theorem {\rm(iv)} for $i=1$. We complete the proof of the main theorem for the case C4.

\subsection{The case C5}

Since the valency of $\Gamma$ is more than $3$, from Lemma \ref{jiben2} (i), we have $k_{1,q-1}=k_{1,q}=2$. Note that $(1,q)$ is mixed. By Theorem \ref{Main1}, $(1,q-1)$ is pure. If $q>2$, from Proposition \ref{Delta(q,q+1)} (i), then $\Gamma\simeq{\rm Cay}(\mathbb{Z}_{4q},\{1,2,2q+1,2q+2\})$. We consider $q=2$ in the following.

By Theorem \ref{Main1} (ii), $p_{(1,2),(1,2)}^{(1,1)}\neq0$. It follows from Lemma \ref{pure-g=2} that $\Gamma_{1,2}\notin\langle\Gamma_{1,1}\rangle$. Suppose $\wz{\partial}(x_{0},x_{1})=(1,2)$ for $x_{0},x_{1}\in V\Gamma$. Then $\partial(F_{2}(x_{0}),F_{2}(x_{1}))=1$ in $\Gamma/\langle\Gamma_{1,1}\rangle$. Since $p_{(1,2),(1,2)}^{(1,1)}\neq0$, we get $\Gamma_{1,1}(x_{0})\cap\Gamma_{1,2}(x_{1})\neq\emptyset$, which implies $\partial(F_{2}(x_{1}),F_{2}(x_{0}))=1$. Hence, $\Gamma/\langle\Gamma_{1,1}\rangle$ is a connected undirected graph. By $k_{1,2}=2$, $\Gamma/\langle\Gamma_{1,1}\rangle\simeq C_{l}$.

Let $(F_{2}(x_{0}),F_{2}(x_{1}),\ldots,F_{2}(x_{l-1}))$ be an undirected circuit. Suppose $l\neq2$. Without loss of generality, we may assume that $(x_{0},x_{1}),(x_{1},x_{2}),(x_{3},x_{2})\in\Gamma_{1,2}$. Then $x_{1}\neq x_{3}$. In view of $\wz{\partial}(x_{0},x_{2})\neq(1,1)$ and Lemma \ref{jiben2} (ii),  one gets $|(\Gamma_{1,2})^{2}|=2$. Since $k_{1,1}=2$, by Lemma \ref{jiben} (i) and (v), we have $p_{(1,2),(1,2)}^{(1,1)}=1$, which implies $x_{3}\in P_{(1,2),(1,2)}(x_{0},x_{2})$ from Lemma \ref{jiben2} (iii). Hence, $\partial(F_{2}(x_{0}),F_{2}(x_{3}))=1$ and $l=4$. Thus, $l=2$ or $4$.

\textbf{Case 1.} $\Gamma/\langle\Gamma_{1,1}\rangle\simeq C_{2}$.

Note that $V\Gamma=F_{2}(x_{0})\dot{\cup}F_{2}(x_{1})$. Let $\tau_{i}$ be an isomorphism from $\textrm{Cay}(\mathbb{Z}_{n},\{1,n-1\})$ to $\Delta_{2}(x_{i})$. Write $\tau_{i}(a)=(a,i)$ for each $a$. Without loss of generality, we may assume $\wz{\partial}_{\Gamma}((0,0),(0,1))=(1,2)$. By Lemma \ref{jiben2} (ii), we get $|(\Gamma_{1,2})^{2}|=1$ or $2$.

\textbf{Case 1.1.} $(\Gamma_{1,2})^{2}=\{\Gamma_{1,1}\}$.

By Lemma \ref{jiben} (i), one has $p_{(1,2),(1,2)}^{(1,1)}=2$, which implies $(1,0),(-1,0)\in\Gamma_{1,2}(0,1)$. It follows from Lemma \ref{pure-g=2} that $\wz{\partial}_{\Gamma}((1,0),(-1,0))=(2,2)$. In view of Lemma \ref{jiben} (ii) and (vi), we get $p_{(1,2),(1,2)}^{(1,1)}p_{(2,1),(1,1)}^{(1,2)}=2+p_{(2,1),(1,2)}^{(2,2)}=4$. By Lemma \ref{jiben} (i) and (v), we obtain $k_{2,2}=1$. It follows from
Lemma \ref{pure-g=2} that $n=4$ and $|V\Gamma|=8$. Since $(\Gamma_{1,2})^{2}=\{\Gamma_{1,1}\}$, by \cite{H}, we obtain $\Gamma\simeq\textrm{Cay}(\mathbb{Z}_{8},\{1,2,5,6\})$.

\textbf{Case 1.2.} $|(\Gamma_{1,2})^{2}|=2$.

Assume that $((0,1),(t,0))\in\Gamma_{1,2}$ and $((0,0),(t,0))\notin\Gamma_{1,1}$. By Theorem \ref{Main1} (iii) and Lemma \ref{pure-g=2}, we have $\wz{\partial}_{\Gamma}((0,0),(t,0))=(2,2)$. Hence, $n>3$. Since $k_{1,1}=2$, we get $p_{(1,2),(1,2)}^{(1,1)}=1$ from Lemma \ref{jiben} (i) and (v). In view of Lemma \ref{jiben2} (iii), we obtain $p_{(1,2),(1,2)}^{(2,2)}=2$ and $k_{2,2}=1$. By Lemma \ref{pure-g=2}, one has $(2,0),(-2,0)\in\Gamma_{2,2}(0,0)$, which implies $n=4$ and $|V\Gamma|=8$. Since $|(\Gamma_{1,2})^{2}|=2$, from \cite{H}, $\Gamma\simeq\textrm{Cay}(\mathbb{Z}_{8},\{1,2,3,6\})$.

\textbf{Case 2.} $\Gamma/\langle\Gamma_{1,1}\rangle\simeq C_{4}$.

Note that $V\Gamma=F_{2}(x_{0})\dot{\cup}F_{2}(x_{1})\dot{\cup}F_{2}(x_{2})\dot{\cup}F_{2}(x_{3})$. Let $\sigma_{i}$ be an isomorphism from $\textrm{Cay}(\mathbb{Z}_{n},\{1,n-1\})$ to $\Delta_{2}(x_{i})$ for each $i$. Write $\tau_{i}(a)=(a,i)$ for any $a$. Without loss of generality, we may assume $\wz{\partial}_{\Gamma}((0,j),(0,j+1))=(1,2)$ for $j=0,1,2$.

Since $(0,j+1)\in P_{(1,2),(1,1)}((0,j),(1,j+1))$, we have $(1,j)$ or $(-1,j)\in\Gamma_{1,2}(1,j+1)$. Without loss of generality, we may assume that $\wz{\partial}_{\Gamma}((1,j),(1,j+1))=(1,2)$. Since $(1,j+1)\in P_{(1,2),(1,1)}((1,j),(2,j+1))$ and $\Gamma/\langle\Gamma_{1,1}\rangle\simeq C_{4}$, one gets $\wz{\partial}_{\Gamma}((2,j),(2,j+1))=(1,2)$. Similarly, $\wz{\partial}_{\Gamma}((a,j),(a,j+1))=(1,2)$ for each $a\in\mathbb{Z}_{n}$ and $j\in\{0,1,2\}$.

By $p_{(1,2),(1,2)}^{(1,1)}\neq0$, we may assume $\wz{\partial}_{\Gamma}((0,1),(1,0))=(1,2)$. Since $(1,0)\in P_{(1,2),(1,1)}((0,1),(2,0))$, we get $(1,1)$ or $(-1,1)\in\Gamma_{2,1}(2,0)$.

\textbf{Case 2.1.} $\wz{\partial}_{\Gamma}((1,1),(2,0))=(1,2)$.

Since $(2,0)\in P_{(1,2),(1,1)}((1,1),(3,0))$ and $\Gamma/\langle\Gamma_{1,1}\rangle\simeq C_{4}$, $\wz{\partial}_{\Gamma}((2,1),(3,0))=(1,2)$. Similarly, $\wz{\partial}_{\Gamma}((a,1),(a+1,0))=(1,2)$ for each $a\in\mathbb{Z}_{n}$. The fact that $p_{(1,2),(1,2)}^{(1,1)}\neq0$ and $\Gamma/\langle\Gamma_{1,1}\rangle\simeq C_{4}$ imply that $(0,2)\in P_{(1,2),(1,2)}((0,1),(-1,1))$. Hence, $((0,0),(0,1),(0,2),(-1,1))$ is a circuit consisting of arcs of type $(1,2)$. In view of Theorem \ref{Main1} (iii), one gets $\wz{\partial}_{\Gamma}((0,0),(0,2))=(2,2)$, which implies $(\Gamma_{1,2})^{2}=\{\Gamma_{1,1},\Gamma_{2,2}\}$ by Lemma \ref{jiben2} (ii). Since $k_{1,1}=2$, from Lemma \ref{jiben} (i) and (v), we obtain $p_{(1,2),(1,2)}^{(1,1)}=1$. In view of Lemma \ref{jiben2} (iii), one has $p_{(1,2),(1,2)}^{(2,2)}=1$ and $k_{2,2}=1$. By Lemma \ref{pure-g=2}, we get $\wz{\partial}_{\Gamma}((0,0),(2,0))=(1,1)$. Since $\Gamma/\langle\Gamma_{1,1}\rangle\simeq C_{4}$, from Theorem \ref{Main1} (i), one obtains $\wz{\partial}_{\Gamma}((0,0),(1,1))=(2,2)$, a contradiction.

\textbf{Case 2.2.} $\wz{\partial}_{\Gamma}((-1,1),(2,0))=(1,2)$.

Since $p_{(1,2),(1,2)}^{(1,1)}\neq0$ and $((-1,0),(-1,2))\notin\Gamma_{1,1}$, we have $\wz{\partial}_{\Gamma}((-1,0),(2,0))=(1,1)$, $n=4$ and $|V\Gamma|=16$. By \cite{H}, $\Gamma\simeq{\rm Cay}(\mathbb{Z}_{4~}\times\mathbb{Z}_{4},\{(0,1),(1,0),(2,0),(0,2)\})$.

We complete the proof of the main theorem for the case C5.

\subsection{The case C6}

By Theorem \ref{Main1} (ii), $p_{(1,q),(1,q)}^{(1,q-1)}\neq0$ and $(1,q-1)$ is pure. In view of Lemma \ref{jiben2} (i), we have $k_{1,1},k_{1,q-1},k_{1,q}\in\{1,2\}$.

\textbf{Case 1.} $k_{1,q-1}=1$.

By Lemma \ref{jiben} (ii), we have $p_{(1,q),(1,q)}^{(1,q-1)}=k_{1,q}$.

\textbf{Case 1.1.} $k_{1,q}=1$.

Since the valency of $\Gamma$ is more than $3$, one has $k_{1,1}=2$. In view of Proposition \ref{Delta(2,q) k 1,q-1=1} and $p_{(1,q),(1,q)}^{(1,q-1)}=1$, $V\Gamma$ has a partition $F_{2,q}(x_{0})\dot{\cup}F_{2,q}(x_{1})$ and there exists an isomorphism $\tau$ from $\textrm{Cay}(\mathbb{Z}_{q}\times\mathbb{Z}_{n},\{(1,0),(0,1),(0,-1)\}$ to $\Delta_{2,q}(x_{0})$ for $n\leq q-(1+(-1)^{q})/2$. Write $\tau(a,b)=(a,b,0)$ for each $(a,b)$. Since $k_{1,q}=1$, $\sigma:F_{2,q}(x_{0})\rightarrow F_{2,q}(x_{1})$, $x\mapsto x'$ is an isomorphism from $\Delta_{2,q}(x_{0})$ to $\Delta_{2,q}(x_{1})$, where $x'\in\Gamma_{1,q}(x)$. Write $\sigma(a,b,0)=(a,b,1)$ for each $(a,b)$. The fact that $p_{(1,q),(1,q)}^{(1,q-1)}=1$ implies $\wz{\partial}_{\Gamma}((a,b,1),(a+1,b,0))=(1,q)$. Thus, $\Gamma$ is isomorphic to one of the digraphs in Theorem \ref{Main} (viii).

\textbf{Case 1.2.} $k_{1,q}=2$.

We claim that $p_{(1,q),(q,1)}^{(1,1)}\neq0$. Let $x,y,z$ be vertices such that $\wz{\partial}(x,y)=(1,1)$ and $\wz{\partial}(y,z)=(1,q-1)$. Since $k_{1,q-1}=1$, by Lemma \ref{pure-g=2}, Theorem \ref{Main1} (i) and Lemma \ref{(1,g-1) (1,1)=(2,g)} (i), one gets $\wz{\partial}(x,z)=(2,q)$. It follows from Lemma \ref{jiben} (i) and Lemma \ref{jiben2} (ii) that $|(\Gamma_{1,q})^{2}|=2$. In view of Lemma \ref{mix-(2,q-2)'} (iii), one gets $p_{(1,q),(1,q)}^{(2,q)}\neq0$, which implies that there exists a vertex $y'\in P_{(1,q),(1,q)}(x,z)$. Since $p_{(1,q),(1,q)}^{(1,q-1)}=2$, we obtain $\wz{\partial}(y,y')=(1,q)$. Thus, our claim is valid.

\textbf{Case 1.2.1.} $k_{1,1}=1$.

Since $k_{1,q-1}=1$, we have $\Gamma_{1,1}\notin\langle\Gamma_{1,q-1}\rangle$. Let $\varphi$ be an isomorphism from $\textrm{Cay}(\mathbb{Z}_{q},\{1\})$ to $\Delta_{q}(x_{0})$ for fixed $x_{0}\in V\Gamma$. Write $\varphi(a)=(a,0,0)$ for any $a\in\mathbb{Z}_{q}$. Since $k_{1,1}=1$, $V\Delta_{2,q}(x_{0})$ has a partition $F_{q}(x_{0})\dot{\cup}F_{q}(x_{1})$. It follows that $\sigma:F_{q}(x_{0})\rightarrow F_{q}(x_{1})$, $x\mapsto x'$ is an isomorphism from $\Delta_{q}(x_{0})$ to $\Delta_{q}(x_{1})$, where $x'\in \Gamma_{1,1}(x)$. Write $\sigma(a,0,0)=(a,1,0)$ for each $a$.

Suppose that there exists $(c,d,0)$ such that $\wz{\partial}_{\Gamma}((0,0,0),(c,d,0))=(1,q)$. The fact that $(1,q-1)$ is pure and $k_{1,q-1}=1$ imply $d=1$ and $c\neq0$. Since $k_{1,1}=1$, by the claim and Lemma \ref{jiben} (v), we get $p_{(1,q),(q,1)}^{(1,1)}=2$, which implies $(c,1,0)\in P_{(1,q),(q,1)}((0,0,0),(0,1,0))$. Then $((0,1,0),(c,1,0),(c+1,1,0),\ldots,(-1,1,0))$ is a circuit of length $q-\hat{c}+1$, a contradiction. Hence, $\Gamma_{1,q}\notin\langle\{\Gamma_{1,1},\Gamma_{1,q-1}\}\rangle$.

Since $p_{(1,q),(1,q)}^{(1,q-1)}=2$, $V\Gamma$ has a partition $F_{2,q}(x_{0})\dot{\cup}F_{2,q}(x_{0}')$. Let $\psi$ be an isomorphism from $\Delta_{2,q}(x_{0})$ to $\Delta_{2,q}(x_{0}')$ such that $\psi(0,0,0)\in\Gamma_{1,q}(0,0,0)$. Write $\psi(a,b,0)=(a,b,1)$ for each $a\in\mathbb{Z}_{q}$ and $b\in\{0,1\}$. Since $p_{(1,q),(1,q)}^{(1,q-1)}=p_{(1,q),(q,1)}^{(1,1)}=2$, we obtain $(a,0,1),(a,1,1)\in\Gamma_{1,q}(a,b,0)$ and $(a+1,0,0),(a+1,1,0)\in\Gamma_{1,q}(a,b,1)$. Then $\Gamma$ is isomorphic to one of the digraphs in Theorem \ref{Main} (v).

\textbf{Case 1.2.1.} $k_{1,1}=2$.

Since $p_{(1,q),(1,q)}^{(1,q-1)}=2$, from Proposition \ref{Delta(2,q) k 1,q-1=1}, $V\Gamma$ has a partition $F_{2,q}(x_{0})\dot{\cup}F_{2,q}(x_{1})$ and there exists an isomorphism $\tau_{i}$ from $\textrm{Cay}(\mathbb{Z}_{q}\times\mathbb{Z}_{n},\{(1,0),(0,1),(0,-1)\})$ to $\Delta_{2,q}(x_{i})$ for $i=0,1$, where $n\leq q-(1+(-1)^{q})/2$. Write $\tau_{i}(a,b)=(a,b,i)$ for each $(a,b)\in\mathbb{Z}_{q}\times\mathbb{Z}_{n}$.

By the claim, we have $p_{(1,q),(q,1)}^{(1,1)}\neq0$. Without loss of generality, we may assume $(0,0,1),(0,-1,1)\in\Gamma_{1,q}(0,0,0)$. In view of $(0,-1,1)\in P_{(1,q),(1,1)}((0,0,0),(0,0,1))$, we may assume $(0,0,1)\in\Gamma_{1,q}(0,1,0)$. Since $(0,0,0)\notin P_{(q,1),(1,q)}((0,0,1),(0,1,1))$ and $p_{(1,q),(q,1)}^{(1,1)}\neq0$, we get $((0,1,0),(0,1,1))\in\Gamma_{1,q}$. Similarly, $(0,b,1),(0,b-1,1)\in\Gamma_{1,q}(0,b,0)$ for each $b$. In view of $p_{(1,q),(1,q)}^{(1,q-1)}=2$. we obtain $(a,b,1),(a,b-1,1)\in\Gamma_{1,q}(a,b,0)$ and $(a+1,b,0),(a+1,b+1,0)\in\Gamma_{1,q}(a,b,1)$ for any $(a,b)\in\mathbb{Z}_{q}\times\mathbb{Z}_{n}$.

Suppose that $c=n/{\rm gcd}(q,n)$ and $c$ is odd. Let $\varphi$ be the mapping from $\Gamma$ to the corresponding digraph in Theorem \ref{Main} (ix) satisfying $\varphi(a,b,i)=(2\hat{a}+i,(2\hat{a}c+ic+i)/2+\hat{b})$. Routinely, $\varphi$ is an isomorphism.

Suppose that $t=q/{\rm gcd}(q,n)$ and $t$ is odd. Let $\psi$ be the mapping from $\Gamma$ to the corresponding digraph in Theorem \ref{Main} (x) such that $\psi(a,b,i)=(2\hat{b}+i,\hat{a}+\hat{b}t+i(1+t)/2)$. Note that $\psi$ is well defined. Assume that $\psi(a,b,i)=\psi(x,y,j)$ for some $(a,b,i)$ and $(x,y,j)$. Since $2\hat{b}+i\equiv2\hat{y}+j~(\textrm{mod}~2n)$, we have $i=j$ and $b=y$. By $\hat{a}+\hat{b}t+i(1+t)/2\equiv\hat{x}+\hat{y}t+j(1+t)/2~(\textrm{mod}~q)$, one gets $a=x$. Therefore, $\psi$ is a bijection. One can verify that
$((x_{1},y_{1},i_{1}),(x_{2},y_{2},i_{2}))$ is an arc if and only if $(\psi(x_{1},y_{1},i_{1}),\psi(x_{2},y_{2},i_{2}))$ is an arc. Hence, $\psi$ is an isomorphism.

\textbf{Case 2.} $k_{1,q-1}=2$.

If $\Gamma_{1,1}\in\Gamma_{1,q-1}\Gamma_{q-1,1}$, by Proposition \ref{Delta(q,q+1)} (i), then $\Gamma\simeq{\rm Cay}(\mathbb{Z}_{4q},\{1,2,2q,2q+1,2q+2\})$ for $q\geq3$. In the following, we consider the case that $\Gamma_{1,1}\notin\Gamma_{1,q-1}\Gamma_{q-1,1}$.

By Proposition \ref{Delta(q,q+1)} (i), there exists an isomorphism $\tau$ from $\textrm{Cay}(\mathbb{Z}_{4q},\{1,2,2q+1,2q+2\})$ to $\Delta_{q,q+1}(x)$ for fixed $x\in V\Gamma$. Write $\tau(a)=(a,0)$ for each $a\in\mathbb{Z}_{4q}$. Observe that $\partial_{\Gamma}((0,0),(b,0))+\partial_{\Gamma}((b,0),(0,0))=q+(1+(-1)^{\hat{b}+1})/2$ for $b\notin\{0,2q\}$. Since $\Gamma_{1,1}\notin\Gamma_{1,q-1}\Gamma_{q-1,1}$, we have $\Gamma_{1,1}\notin\langle\{\Gamma_{1,q-1},\Gamma_{1,q}\}\rangle$.

\textbf{Case 2.1.} $k_{1,1}=1$.

Observe that $V\Gamma$ has a partition $F_{q,q+1}(x)\dot{\cup}F_{q,q+1}(x')$. Note that $\sigma:F_{q,q+1}(x)\rightarrow F_{q,q+1}(x')$, $y\mapsto y'$ is an isomorphism from $\Delta_{q,q+1}(x)$ to $\Delta_{q,q+1}(x')$, where $y'\in\Gamma_{1,1}(y)$. Write $\sigma(a,0)=(a,1)$ for each $a$. If $q=3$, then $(6,0),(3,1),(9,1)\in\Gamma_{3,3}(0,0)$, a contradiction. Hence, $\Gamma$ is isomorphic to one of the digraphs in Theorem \ref{Main} (vi) for $i=0$.

\textbf{Case 2.2.} $k_{1,1}=2$.

By Proposition \ref{pure} and Lemma \ref{Gamma 2,q} (i), (iv), one gets $A_{1,q-1}A_{1,1}=2A_{2,q}$ and $A_{1,q}A_{1,1}=2A_{2,q+1}$. Hence, $V\Gamma$ has a partition $F_{q,q+1}(x)\dot{\cup}F_{q,q+1}(x')$. Let $\varphi$ be an isomorphism from $\Delta_{q,q+1}(x)$ to $\Delta_{q,q+1}(x')$ such that $\varphi(0,0)\in\Gamma_{1,1}(0,0)$. Write $\varphi(a,0)=(a,1)$ for each $a$. Since $A_{1,q-1}A_{1,1}=2A_{2,q}$ and $A_{1,q}A_{1,1}=2A_{2,q+1}$, we have $(a,1),(a+2q,1)\in\Gamma_{1,1}(a,0)$ for each $a$. If $2<q<5$, then $(2q,0),(q,0),(3q,0)\in\Gamma_{2,2}(0,0)$, contrary to Lemma \ref{jiben2} (i). Therefore, $\Gamma$ is isomorphic to one of the digraphs in Theorem \ref{Main} (vi) for $i=1$.

We complete the proof of the main theorem for the case C6.

\section*{Acknowledgement}
 Y. Yang is supported by the Fundamental Research Funds for the Central Universities (Grant No. 2652017141), B. Lv is supported by NSFC (11501036), K. Wang is supported by NSFC (11671043).

\newpage

\begin{table}
\renewcommand\arraystretch{1.15}
\caption{\quad Two way distance of digraphs in Theorem \ref{Main}} \vspace{0.2cm}

\setlength{\tabcolsep}{2pt}
\begin{tabular}{|c|c|c|}
\hline
$\Gamma$ & Conditions & $\wz{\partial}((0,0),(a,b))$ with $(a,b)\neq(0,0)$ \\\hline
\multirow{2}{*}{(iv)} & $a\neq0$ & $(\beta(\hat{b})+\hat{a},q+\beta(\hat{b})-\hat{a})$ \\\cline{2-3}
& $a=0$ & $(\lceil\frac{\hat{b}}{2}\rceil+(-1)^{\hat{b}}\lceil\frac{\hat{b}-1}{2}\rceil i,\lceil2-\frac{\hat{b}}{2}\rceil+(-1)^{\hat{b}}\lceil\frac{3-\hat{b}}{2}\rceil i)$ \\\hline
\multirow{3}{*}{(v)}& $2\nmid\hat{a}$ & $(\frac{\hat{a}+1}{2},q-\frac{\hat{a}-1}{2})$ \\\cline{2-3}
& $(a,b)\neq(0,1)$ and $2\mid\hat{a}$ & $(\hat{b}+\frac{\hat{a}}{2},q+\hat{b}-\frac{\hat{a}}{2})$ \\\cline{2-3}
& $(a,b)=(0,1)$ & $(1,1)$ \\\hline
\multirow{4}{*}{(vi)} & $0<\hat{a}<2q$ & $(\frac{\hat{a}+2\hat{b}+\beta(\hat{a})}{2},q-\frac{\hat{a}-2\hat{b}-\beta(\hat{a})}{2})$ \\\cline{2-3}
& $\hat{a}>2q$ & $(\frac{\hat{a}+2\hat{b}+\beta(\hat{a})}{2}-q,2q-\frac{\hat{a}-2\hat{b}-\beta(\hat{a})}{2})$ \\\cline{2-3}
& $a=2q$ & $(q^{1-i}+\hat{b}+(-1)^{\hat{b}}i,q^{1-i}+\hat{b}+(-1)^{\hat{b}}i)$ \\\cline{2-3}
& $(a,b)=(0,1)$ & $(1,1)$ \\\hline
\multirow{2}{*}{(vii)} & $a\neq0$ & $(\beta(\hat{b})+\frac{\hat{a}+\beta(\hat{a})}{2},q+\beta(\hat{b})-\frac{\hat{a}-\beta(\hat{a})}{2})$ \\\cline{2-3}
& $a=0$ & $(\lceil\frac{\hat{b}}{2}\rceil+\lceil\frac{\hat{b}-1}{2}\rceil i,\lceil2-\frac{\hat{b}}{2}\rceil+\lceil\frac{3-\hat{b}}{2}\rceil i)$ \\\hline
\multirow{4}{*}{(viii)}& $a=0$ and $\hat{b}\leq\frac{n}{2}$ & $(\hat{b},\hat{b})$ \\\cline{2-3}
& $a=0$ and $\hat{b}>\frac{n}{2}$ & $(n-\hat{b},n-\hat{b})$ \\\cline{2-3}
& $a\neq0$ and $\hat{b}\leq\frac{n}{2}$ & $(\hat{b}+\frac{\hat{a}+\beta(\hat{a})}{2},\hat{b}+q-\frac{\hat{a}-\beta(\hat{a})}{2})$ \\\cline{2-3}
& $a\neq0$ and $\hat{b}>\frac{n}{2}$ & $(n-\hat{b}+\frac{\hat{a}+\beta(\hat{a})}{2},n-\hat{b}+q-\frac{\hat{a}-\beta(\hat{a})}{2})$ \\\hline
\multirow{4}{*}{(ix)}& $a=0$ and $v_{a,b}\leq\frac{n}{2}$  & $(v_{a,b},v_{a,b})$ \\\cline{2-3}
& $a=0$ and $v_{a,b}>\frac{n}{2}$ & $(n-v_{a,b},n-v_{a,b})$ \\\cline{2-3}
& $a\neq0$ and $v_{a,b}\leq\frac{n-\beta(\hat{a})}{2}$ & $(v_{a,b}+\frac{\hat{a}+\beta(\hat{a})}{2},v_{a,b}+q-\frac{\hat{a}-\beta(\hat{a})}{2})$ \\\cline{2-3}
& $a\neq0$ and $v_{a,b}>\frac{n-\beta(\hat{a})}{2}$ & $(n-v_{a,b}+\frac{\hat{a}-\beta(\hat{a})}{2},n-v_{a,b}+q-\frac{\hat{a}+\beta(\hat{a})}{2})$ \\\hline
\multirow{4}{*}{(x)}& $u_{a,b}=0$ and $v_{a}\leq\frac{n}{2}$  & $(v_{a},v_{a})$ \\\cline{2-3}
& $u_{a,b}=0$ and $v_{a}>\frac{n}{2}$ & $(n-v_{a},n-v_{a})$ \\\cline{2-3}
& $u_{a,b}\neq0$ and $v_{a}\leq\frac{n-\beta(u_{a,b})}{2}$ & $(v_{a}+\frac{u_{a,b}+\beta(u_{a,b})}{2},v_{a}+q-\frac{u_{a,b}-\beta(u_{a,b})}{2})$ \\\cline{2-3}
& $u_{a,b}\neq0$ and $v_{a}>\frac{n-\beta(u_{a,b})}{2}$ & $(n-v_{a}+\frac{u_{a,b}-\beta(u_{a,b})}{2},n-v_{a}+q-\frac{u_{a,b}+\beta(u_{a,b})}{2})$ \\\hline
\end{tabular}
\begin{center}
$\beta(q)=(1+(-1)^{q+1})/2$, \ $v_{a}=(\hat{a}-\beta(\hat{a}))/2$,\\
$0\leq v_{a,b}<n$ and $v_{a,b}\equiv\hat{b}-(\hat{a}c+\beta(\hat{a}))/2~({\rm mod}~n)$,\\
$0\leq u_{a,b}<q$ and $u_{a,b}\equiv2\hat{b}-\beta(\hat{a})t-2tv_{a}~({\rm mod}~q)$.
\end{center}
\end{table}
\end{CJK*}

\end{document}